\documentclass[a4paper,12pt]{amsart}
\usepackage{amscd}
\usepackage{amsfonts}
\usepackage{amsmath}
\usepackage{amssymb}
\usepackage{latexsym}
\usepackage{amsthm,color}
\usepackage[all]{xy}
\usepackage{epsfig}
\usepackage{graphicx}
\usepackage{color}
\usepackage{tikz}
\usetikzlibrary{calc}
\newtheorem{proposition}{Proposition}
\newtheorem{lemma}{Lemma}[section]

\newtheorem{theorem}{Theorem}
\newtheorem{corollary}{Corollary}

\theoremstyle{definition}
\newtheorem{definition}{Definition}
\theoremstyle{definition}
\newtheorem{example}{Example}%[section]
\theoremstyle{definition}

\theoremstyle{definition}

\newcommand{\R}{\mathbb{R}}

\begin{document}
\title[
Anti-orthotomics of frontals and their applications]
{Anti-orthotomics of frontals \\ 
and their applications}
%%%%%%%%%%%%%%%%%%%%%%%%%%%%%%%%%%%%%%%%%%%%%%%%%%%%%%%%%%%%%%%%%%% 
%%%%%%%%%%%%%%%%%%%%%%%%%%%%%%%%%%%%%%%%%%%%%%%%%%%%%%%%%%%%%%%%%%%
%\thanks{\color{black} This work is partially supported 
%by JSPS and CAPES 
%under the Japan--Brazil research cooperative program and 
%JSPS KAKENHI Grant Number 26610035.}
%%%%%%%%%%%%%%%%%%%%%%%%%%%%%%%%%%%%%%%%%%%%%%  
\author[S.~Janeczko]{Stanis{\l}aw Janeczko}
\address{Institute of Mathematics, 
Polish Academy of Sciences, 
ul. Sniadeckich 8, 00-956 Warsaw, POLAND \\ 
and \\  
Faculty of Mathematics and Information Science, 
Warsaw University of Technology, 
ul. Koszykowa 75, 00-662 Warsaw, POLAND
}
\email{S.Janeczko@mini.pw.edu.pl}
\author[T.~Nishimura]{Takashi Nishimura
%\thanks{Corresponding author.} 
}
\address{
Research Institute of Environment and Information Sciences,  
Yokohama National University, 
240-8501 Yokohama, JAPAN}
\email{nishimura-takashi-yx@ynu.ac.jp}
%%%%%%%%%%%%%%%%%%%%%%%%%%%%%%%%%%%%%%%%%%%%%%%  
%%%%%%%%%%%%%%%%%%%%%%%%%%%%%%%%%%%%%%%%%%%%%%% 
\begin{abstract}
Let $f: N^n\to \mathbb{R}^{n+1}$ be 
a frontal with its Gauss mapping  
$\nu: N\to S^n$ and let 
$P\in \mathbb{R}^{n+1}$ be a point such that 
$(f(x)-P)\cdot \nu(x) \ne 0$ for any $x\in N$.    
In this paper, for   
the mapping $\widetilde{f}: N\to \mathbb{R}^{n+1}$ 
defined by 
\[
\widetilde{f}(x)=f(x)-\frac{||f(x)-P||^2}{2(f(x)-P)
\cdot \nu(x)}\nu(x),  
\]
the following four are shown.   
%\begin{enumerate}
(1) $\widetilde{f}$ is 
a frontal  with its Gauss mapping   
$\widetilde{\nu}(x)=\frac{f(x)-P}{||f(x)-P||}$ 
at $\widetilde{f}(x)$.   
(2) $\widetilde{f}$ 
is the unique anti-orthotomic of $f$ relative to $P$.   
(3) The property 
$(\widetilde{f}(x)-P)\cdot \widetilde{\nu}(x)\ne 0
$ holds for any $x\in N$.   
(4) The equality 
$||\widetilde{f}(x)-P||=||\widetilde{f}(x)-f(x)||$ 
holds for any $x\in N$.  
%\end{enumerate}
\par 
Moreover, three applications of the main result are given.  
As the first application, a generalization of Cahn-Hoffman 
vector formula is given.   
The second application is to clarify  an optical meaning 
of anti-orthotomics.   The third application gives a criterion 
to be a front for a given frontal.   
% , a detailed investigation 
%is given for iterated reflections 
%by the anti-orthotomic mirror of the unit sphere 
%$S^n$ relative to a point $P$ $(||P||<1)$.    
\end{abstract}
%, Gauss mapping,  
%Immersion, Hedgehog, 
%Iterated reflections, 
%Mirror, Light-source, Iterated reflection %sequence, 
%Periodic point, Attractor of iterated reflection sequence, }
%\thanks{}

\date{}

\maketitle

%\maketitle

\section{Introduction\label{section 1}}
Throughout this paper, let $n$, $N$ be 
a positive integer and an 
$n$-dimensional $C^\infty$ manifold 
without boundary respectively.    
Moreover, all mappings in this paper are 
of class $C^\infty$ unless otherwise 
stated.   
\par  
A mapping $f: N\to \mathbb{R}^{n+1}$ is called a 
{\it frontal} if there exists a mapping 
$\nu: N\to \mathbb{R}^{n+1}$ such that 
the following two conditions are 
satisfied, where the dot in the center stands 
for the scalar product of 
two vectors in $\mathbb{R}^{n+1}$ and two vector spaces 
$T_{f(x)}\mathbb{R}^{n+1}$ and $\mathbb{R}^{n+1}$ 
are identified:   
\begin{enumerate}
\item[(1)] $\nu(x)\cdot \nu(x)=1$, i.e. $\nu(x)\in S^n$ 
for any $x\in N$.  
\item[(2)] $df_x({\bf v})\cdot \nu(x)=0$ for any $x\in N$ and 
any ${\bf v}\in T_xN$.    
\end{enumerate}  
By the above conditions $(1)$ and $(2)$, it is natural to call 
$\nu: N\to S^n$ the \textit{Gauss mapping} of $f$.   
In this paper, sometimes, even the mapping 
$(f, \nu): N\to \mathbb{R}^{n+1}\times S^n$ is called a 
frontal.   
The notion of frontal was independently introduced in several  
literature (e.g. \cite{fujimorisajiumeharayamada, 
ishikawa2000, sajiumeharayamada, zakalyukinlurbatskii}) 
and it has been rapidly 
and intensively  investigated 
(%for instance, 
see %excellent survey articles 
\cite{ishikawa}).    
\begin{definition}\label{orthotomic}
Let $\widetilde{f}: N\to \mathbb{R}^{n+1}$ be 
a frontal with its Gauss mapping $\widetilde{\nu}: 
N\to S^n$ and  
%satisfying conditions of frontal for ${f}$.   
let $P$ be a point of $\mathbb{R}^{n+1}$.    
\begin{enumerate}
\item[(1)] A mapping ${f}: N\to \mathbb{R}^{n+1}$ is called the  
\textit{orthotomic} of $\widetilde{f}$ relative to $P$ 
if the following equality holds 
for any $x\in N$.     
\[
{f}(x)=2\left(\left(\widetilde{f}(x)-P\right)\cdot 
\widetilde\nu(x)\right)\widetilde\nu(x)+P.   
\] 
\item[(2)] A mapping ${g}: N\to \mathbb{R}^{n+1}$ is called the  
\textit{pedal} of $\widetilde{f}$ relative to $P$ 
if the following equality holds for any $x\in N$.     
\[
{g}(x)=\left(\left(\widetilde{f}(x)-P\right)\cdot 
\widetilde\nu(x)\right)\widetilde\nu(x)+P.   
\] 
\end{enumerate}
\end{definition}
%%%%%%%%%%%%%%%%%%%%%%%%%%%%%%%%%%%%%%%%%%%%%%   
\begin{figure}
\begin{center}
\includegraphics[width=6cm]
{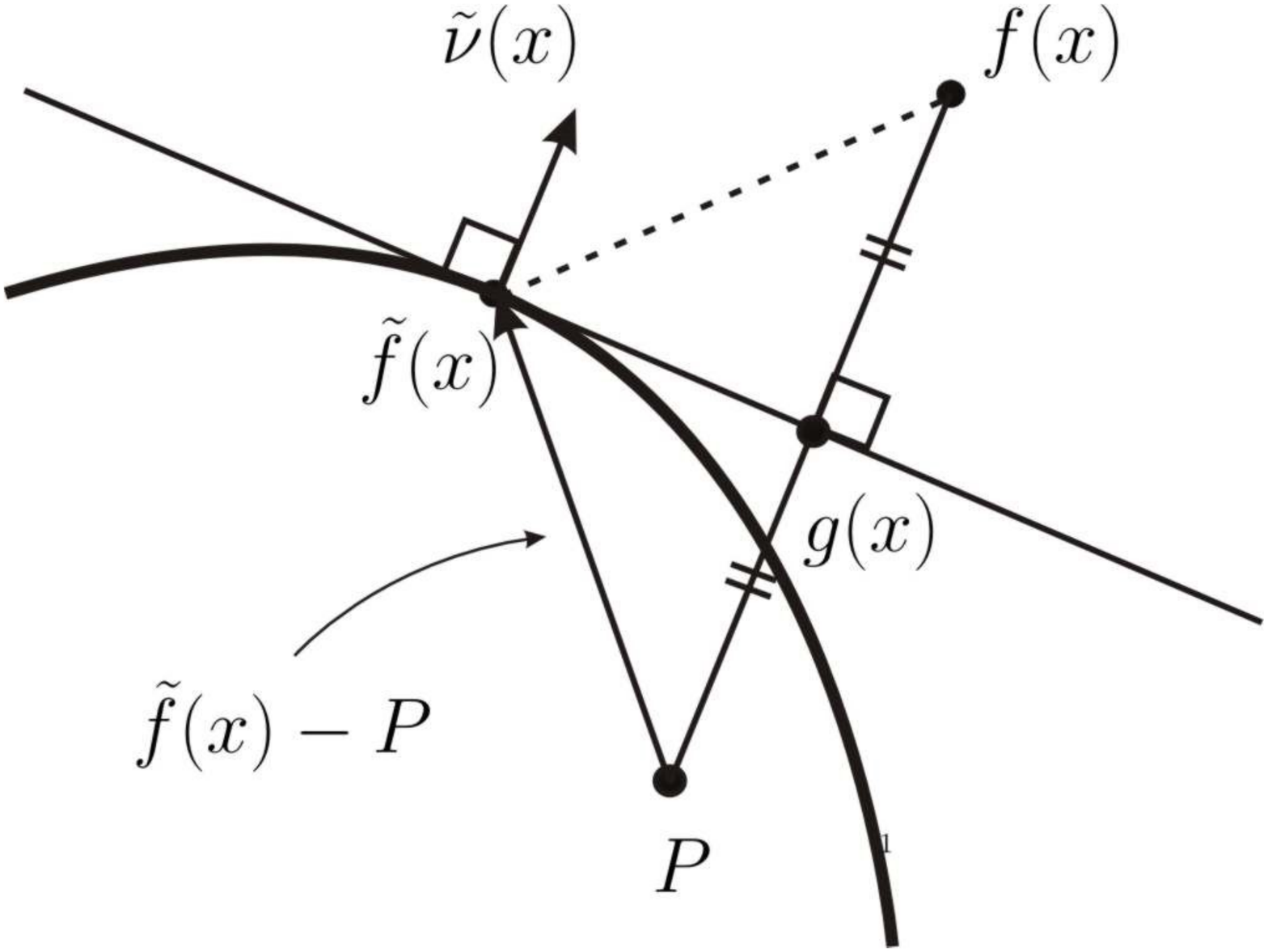}
\caption{%Left:
Orthotomic $f$ and pedal $g$ 
of $(\widetilde{f}, \widetilde{\nu})$ 
relative to $P$. 
% Right: 
}
\label{figure 1}
\end{center}
\end{figure}
%%%%%%%%%%%%%%%%%%%%%%%%%%%%%%%%%%%%%%%%%%%%%%% 
%%%%%%%%%%%%%%%%%%%%%%%%%%%%%%%%%%%%%%% 
%%%%%%%%%%%%%%%%%%%%%%%%%%%%%%%%%%%%%%% 
%\vspace{7cm} 
%%%%%%%%%%%%%%%%%%%%%%%%%%%%%%%%%%%%%%% 
%%%%%%%%%%%%%%%%%%%%%%%%%%%%%%%%%%%%%%% 
{\color{black}Figure \ref{figure 1} clearly 
illustrates the relation between the orthotomic $f$ and 
the pedal $g$ of $\widetilde{f}$ relative to $P$}.     
%%%%%%%%%%%%%%%%%%%%%%%%%%%%%%%%%%%%%%% 
The following Proposition \ref{proposition 1} guarantees 
that the orthotomic of a given 
frontal %relative to a given generic point 
is a frontal.  
\begin{proposition}\label{proposition 1} 
Let $(\widetilde{f}, \widetilde{\nu}): 
N\to \mathbb{R}^{n+1}$ be a frontal 
%with unit normal vector $\nu(x)$ 
and let $P$ be a point of $\mathbb{R}^{n+1}$ such that 
the following condition is satisfied for any $x\in N$.   
\[
\left(\widetilde{f}(x)-P\right)\cdot \widetilde{\nu}(x)\ne 0. 
%P\not\in f(N) 
\]    
Then, the orthotomic of 
$\widetilde{f}$ relative to $P$ defined by 
\[
{f}(x)=2\left((\widetilde{f}(x)-P)\cdot 
\widetilde{\nu}(x)\right)\widetilde{\nu}(x)+P
\] 
is a frontal 
with its Gauss mapping  
${\nu}(x)=
\frac{f(x)-\widetilde{f}(x)}{||f(x)-\widetilde{f}(x)||}$. 
Moreover, the condition  
\[
\left({f}(x)-P\right)\cdot {\nu}(x)\ne 0
\] 
holds for any $x\in N$.  
\end{proposition} 
%In the case $n=1$, Proposition \ref{proposition 1} 
%has been shown 
%in \cite{lipei}.    
\noindent 
Proposition \ref{proposition 1} will be proved in Section 2.   
By definition, it is clear that 
for the pedal $g$, 
$f=2g-P$ is the orthotomic.    
Thus, it is clear that Proposition \ref{proposition 1} yields 
the following corollary which is a generalization of 
\cite{lipei}.   
\begin{corollary}\label{corollary 1} 
Let $(\widetilde{f}, \widetilde{\nu}): 
N\to \mathbb{R}^{n+1}$ be a frontal 
%with unit normal vector $\nu(x)$ 
and let $P$ be a point of $\mathbb{R}^{n+1}$ such that 
the following condition is satisfied for any $x\in N$.   
\[
\left(\widetilde{f}(x)-P\right)\cdot \widetilde{\nu}(x)\ne 0. 
%P\not\in f(N) 
\]    
Then, the pedal of 
$\widetilde{f}$ relative to $P$ defined by 
\[
{g}(x)=\left((\widetilde{f}(x)-P)\cdot 
\widetilde{\nu}(x)\right)\widetilde{\nu}(x)+P
\] 
is a frontal 
with its Gauss mapping  
${\nu}(x)=
\frac{2g(x)-P-\widetilde{f}(x)}
{||2g(x)-P-\widetilde{f}(x)||}$. 
Moreover, the condition  
\[
\left({g}(x)-P\right)\cdot {\nu}(x)\ne 0
\] 
holds for any $x\in N$.  
\end{corollary} 
Notice that in the case that $\widetilde{f}$ 
is a plane regular curve, it is well-known that 
$\widetilde{f}(x)-f(x)$ 
is a normal vector to ${f}$  at ${f}(x)$ 
(for instance, see \cite{brucegiblin}).  
Therefore, 
a part of Proposition \ref{proposition 1} may be regarded 
as just a generalization of the classical result 
to frontals of general dimension.    
\par 
Notice also that 
even if $\widetilde{f}: N\to \mathbb{R}^{n+1}$ 
is non-singular, 
the condition 
\lq\lq $(\widetilde{f}(x)-P)\cdot \widetilde{\nu}(x)\ne 0$ 
for any $x\in N$\rq\rq\;  
seems not so mild.     In other words, 
even when $\widetilde{f}: N\to \mathbb{R}^{n+1}$ 
is an embedding, 
if the image of Gaussian curvature function of 
$\widetilde{f}(N)$ 
is a large interval 
containing zero as an interior point, 
then there are no points $P$ 
satisfying the condition 
\lq\lq $(\widetilde{f}(x)-P)\cdot \widetilde{\nu}(x)\ne 0$ 
for any $x\in N$\rq\rq.           
On the other hand, 
if $\widetilde{f}: N\to \mathbb{R}^{n+1}$ is an embedding and  
the Gaussian curvature of $\widetilde{f}(N)$ 
is always positive, then 
the set 
$\{P\in \R^{n+1}\; |\;  (\widetilde{f}(x)-P)\cdot 
\widetilde{\nu}(x)\ne 0 
\mbox{ for any }x\in N\}$ 
is a non-empty open set.     
Moreover, 
the assumption  \lq\lq $\widetilde{f}: N\to \mathbb{R}^{n+1}$ 
is an embedding and  
the Gaussian curvature of $\widetilde{f}(N)$ 
is always positive\rq\rq\; 
seems to be 
common for the study of orthotomics and pedals.   
(for instance, see \cite{alamocriado, brucegiblingibson}).     
Therefore, the assumption given in Proposition 
\ref{proposition 1} generalizes 
the common assumption for the study of orthotomics 
and pedals.   
\par
\begin{figure}
\begin{center}
\includegraphics[width=5.8cm]
{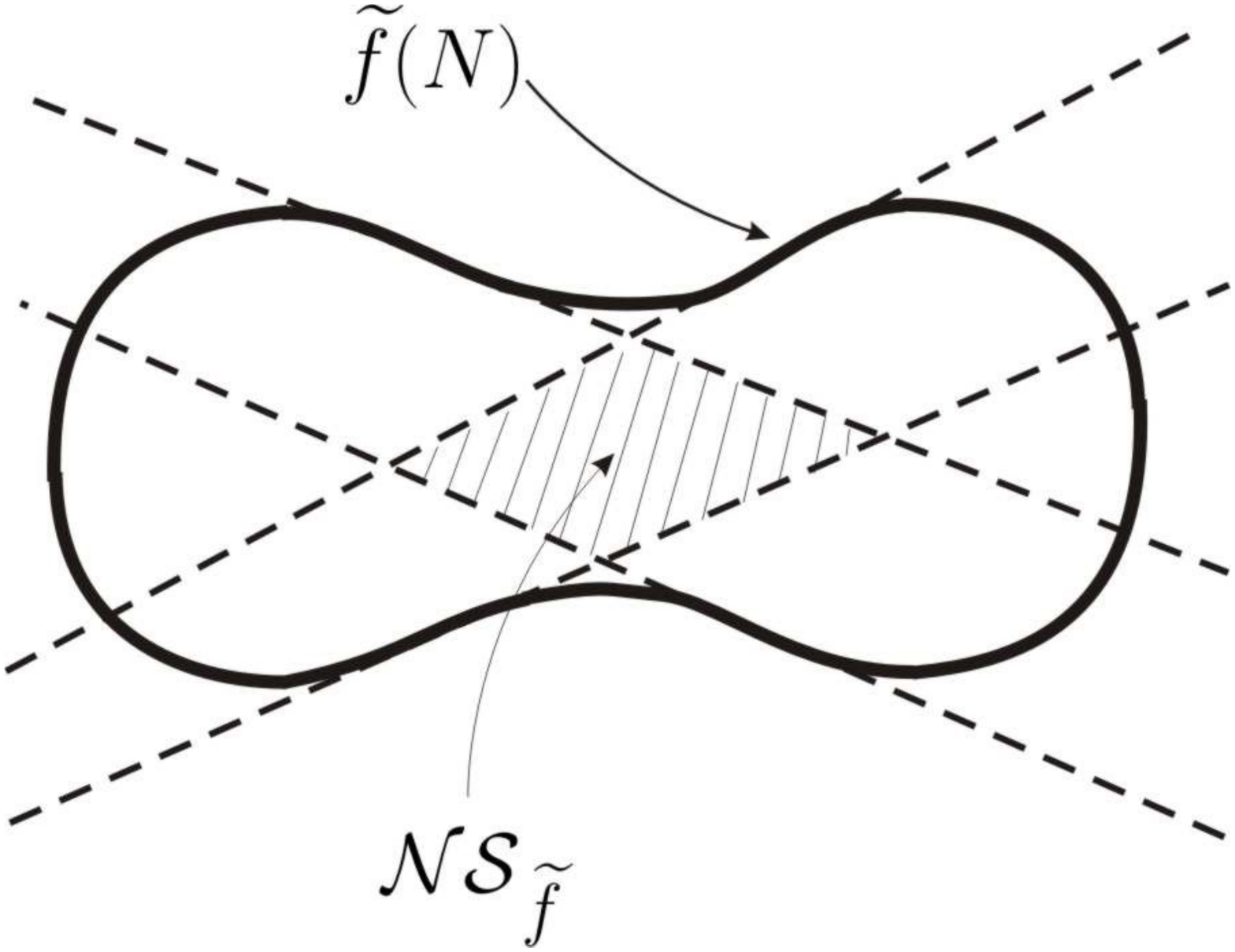}
\qquad 
\includegraphics[width=5.6cm]
{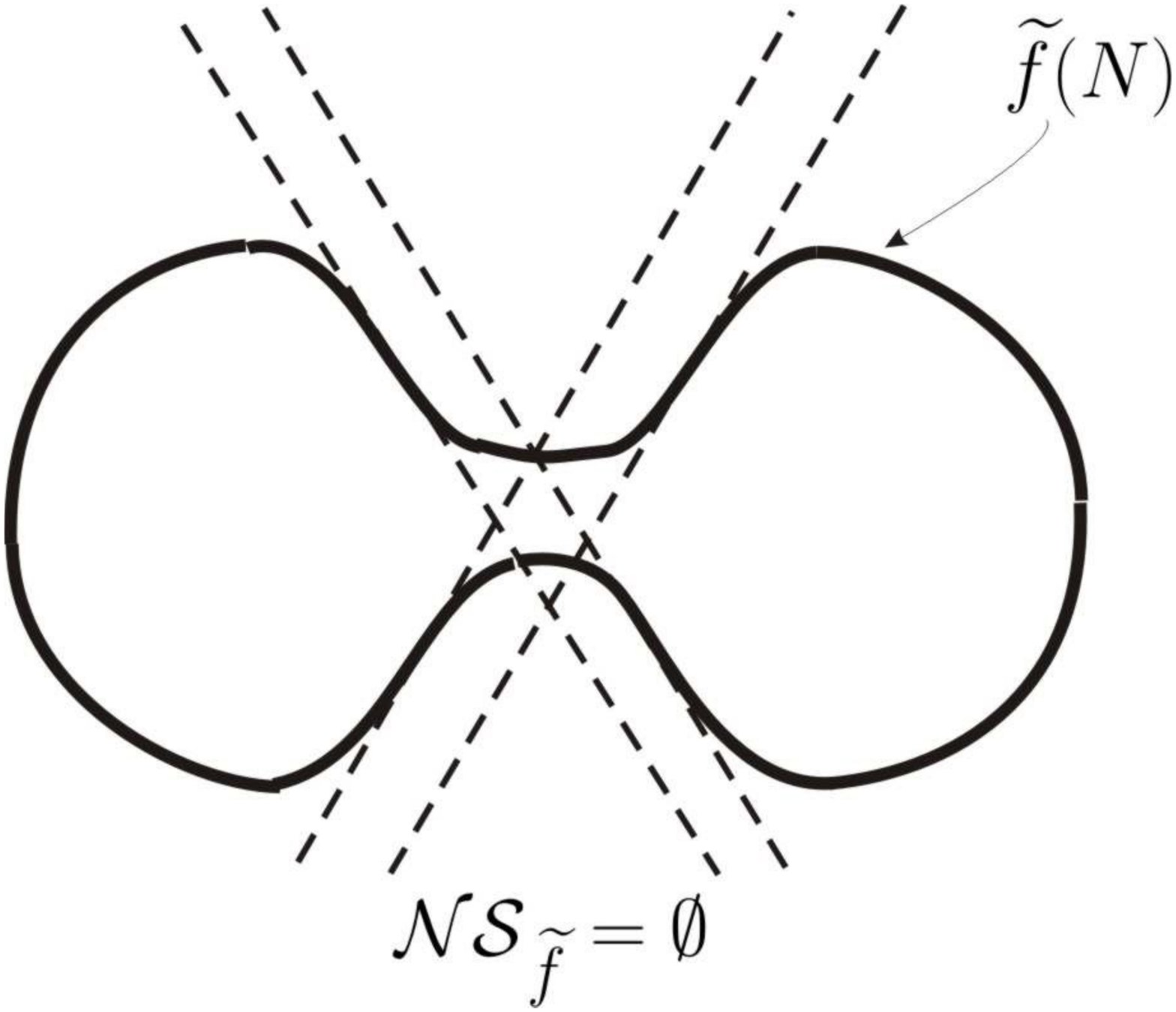}
\caption{
Left: 
$\mathcal{NS}_{\widetilde{f}}$ is not empty.    
Right: $\mathcal{NS}_{\widetilde{f}}$ is empty.
} 
\label{figure 2}
\end{center}
\end{figure}

 
The same condition as the assumption 
of Proposition \ref{proposition 1} 
has been already introduced 
%%%%%%%%%%%%%%%%%%%%%%%%%%%%%%%%%
{\color{black} by J.~W. Bruce and P.~J.~Giblin in 
\cite{brucegiblin} 7.14 in the case that 
$\widetilde{f}: I\to \mathbb{R}^2$ is regular; 
and also 
} 
%%%%%%%%%%%%%%%%%%%%%%%%%%%%%%%%%
by the 
second author (\cite{nishimurasakemi} in the case that 
$\widetilde{f}: N\to S^{n+1}$ is an immersion, 
\cite{nishimurasakemi2} in the case that 
$\widetilde{f}: S^n\to \mathbb{R}^{n+1}$ 
is a Legendrian map and 
\cite{kagatsumenishimura} in the case that $n=1$ 
and $\widetilde{f}: S^1\to \mathbb{R}^2$ is 
an embedding).    
Namely, in \cite{kagatsumenishimura} 
%(resp., \cite{nishimurasakemi}) 
the following set, called {\it no-silhouette} of $\widetilde{f}$ 
and denoted by 
$\mathcal{NS}_{\widetilde{f}}$, 
is defined.   
% and the mapping 
%$\widetilde{f}: S^1\to \mathbb{R}^2$ 
%with non-empty $\mathcal{NS}_{\widetilde{f}}$ 
%have been studied. 
\[
\mathcal{NS}_{\widetilde{f}}
=\left\{P\in \mathbb{R}^{2}\; \left| \; 
\mathbb{R}^{2}- \bigcup_{x\in S^1}
\left(\widetilde{f}(x)+
d\widetilde{f}_x(T_xS^1)\right)\right.\right\}.   
\]  
For a frontal $\widetilde{f}: N\to \mathbb{R}^{n+1}$ 
with its Gauss mapping $\widetilde{\nu}$, 
the notion of no-silhouette $\mathcal{NS}_{\widetilde{f}}$ 
can be naturally generalized as follows.    
{\color{black}The optical meaning of no-silhouette 
is illustrated by Figure \ref{figure 2}.}         
\[
\mathcal{NS}_{\widetilde{f}}
=\left\{\left.P\in \mathbb{R}^{n+1}\; \right| \; 
\left(\widetilde{f}(x)-P\right)\cdot 
\widetilde{\nu}(x)\ne 0\mbox{ for any }x\in N  
\right\}.   
\]   
%By using the notion of no-silhouette $\mathcal{NS}_f$, 
%the condition \lq\lq $(f(x)-P)\cdot \nu(x)\ne 0$ for any 
%$x\in N$\rq\rq can be simply expressed 
%as $P\in \mathcal{NS}_f$.   
%\par 
%\medskip 
%%%%%%%%%%%%%%%%%%%%%%%%%%%%%%%%%%%%%%%%%%%%%%   
%%%%%%%%%%%%%%%%%%%%%%%%%%%%%%%%%%%%%%%%%%%%%%% 
%%%%%%%%%%%%%%%%%%%%%%%%%%%%%%%%%%%%%%  
%%%%%%%%%%%%%%%%%%%%%%%%%%%%%%%%%%%%%% 
%\vspace{7cm}
%%%%%%%%%%%%%%%%%%%%%%%%%%%%%%%%%%%%%% 
%%%%%%%%%%%%%%%%%%%%%%%%%%%%%%%%%%%%%% 
\begin{definition}\label{anti-orthotomic}
\begin{enumerate}
\item[(1)] 
Let ${f}: N\to \mathbb{R}^{n+1}$ be a $C^\infty$ 
mapping    
%frontal with its Gauss mapping ${\nu}: 
%N\to \mathbb{R}^{n+1}$ and  
%satisfying conditions of frontal for ${f}$.   
and let $P$ be a point of $\mathbb{R}^{n+1}$.    
A frontal  
$\widetilde{f}: N\to \mathbb{R}^{n+1}$ 
with its Gauss mapping  
$\widetilde{\nu}: N\to S^n$ 
is called the \textit{anti-orthotomic} of $f$ relative to $P$ 
if the following equality holds 
for any $x\in N$.     
\[
{f}(x)=2\left(\left(\widetilde{f}(x)-P\right)\cdot 
\widetilde\nu(x)\right)\widetilde\nu(x)+P.   
\] 
\item[(2)] 
Let $g: N\to \mathbb{R}^{n+1}$ be a $C^\infty$ 
mapping    
%frontal with its Gauss mapping ${\nu}: 
%N\to \mathbb{R}^{n+1}$ and  
%satisfying conditions of frontal for ${f}$.   
and let $P$ be a point of $\mathbb{R}^{n+1}$.    
A frontal  
$\widetilde{f}: N\to \mathbb{R}^{n+1}$ 
with its Gauss mapping  
$\widetilde{\nu}: N\to S^n$ 
is called the \textit{negative pedal} of $g$ relative to $P$ 
if the following equality holds 
for any $x\in N$.     
\[
g(x)=\left(\left(\widetilde{f}(x)-P\right)\cdot 
\widetilde\nu(x)\right)\widetilde\nu(x)+P.   
\] 
\end{enumerate}
\end{definition} 
By definition, if $f=2g-P$, then the anti-orthotomic 
of $f$ relative to $P$ is exactly the same as 
the negative pedal of $g$ relative to $P$.     
Depending on situations, sometimes,  
the negative pedal is also called the 
\textit{primitive} 
(for example, see \cite{arnold}) or the 
\textit{Cahn-Hoffman map} 
(for instance, see \cite{hondakoisotanaka, 
jikumatsukoiso, koiso}).     
In Geometric Optics, the notions of anti-orthotomic is  
very important (for example, see 
\cite{alamocriado, brucegiblin, brucegiblingibson}), and 
for the study of Wulff shape, the notion of negative pedal 
is the core concept (for instance, see 
\cite{hannishimura, jikumatsukoiso, koiso, morgan, taylor}).    
\par 
By definition, it is clear that an anti-orthotomic 
(resp., a negative pedal)  
is a solution frontal    
for a given orthotomic equation (resp., pedal equation).    
Therefore, obtaining anti-orthotomics or negative 
pedals may be considered as inverse problems.   
It seems that, except for the case that the Gauss 
mapping $\widetilde{\nu}$ of $\widetilde{f}$ 
is non-singular (i.e. the case that 
$\widetilde{f}: N\to \mathbb{R}^{n+1}$ 
is an embedding and  
the Gaussian curvature of $\widetilde{f}(N)$ 
is always positive), 
such inverse problems 
have been usually investigated  
only by solving simultaneous function equations 
for the envelopes.       
\begin{example}\label{example 1}
Let $g: \mathbb{R}\to \mathbb{R}^2$ be the mapping 
defined by 
$g(\theta)=\left(\cos\left(\theta^3\right), 
\sin\left(\theta^3\right)\right)$. 
Define $\nu: \mathbb{R}\to \mathbb{R}^2$ by $\nu=g$.    
Then, $(g, \nu): \mathbb{R}\to \mathbb{R}^2\times S^1$ 
is a frontal.    
Set $P=(0,0)$ and $\widetilde{f}=\widetilde{\nu}=g$.     
Then, clearly, 
$(\widetilde{f}, \widetilde{\nu}): 
\mathbb{R}\to \mathbb{R}^2\times S^1$  
is the negative pedal of $g$ relative to $P$ and 
$\widetilde{f}(\mathbb{R})=S^1$.    
\par 
On the other hand, the function  
$\Phi: \mathbb{R}^2\times \mathbb{R}\to \mathbb{R}$ 
given by 
\[
\Phi(X, Y, \theta)=\nu(\theta)\cdot ((X,Y)-g(\theta))
\]       
defines the one-parameter family of affine 
tangent lines to $g(\mathbb{R})=S^1$.     
And the solution figure of the simultaneous equation 
$\Phi=0, \frac{\partial\Phi}{\partial\theta}=0$ 
is $S^1\cup \{(1, Y)\; |\; Y\in \mathbb{R}\}$.    
\end{example}
\begin{example}\label{example 2}
Let $x: \mathbb{R}\to \mathbb{R}$ be a $C^\infty$ 
periodic function of period $8$ satisfying 
the following condition for any $n$. 
\[
\left\{
\begin{array}{cc}
x(t)=1 & (\mbox{if }\; 8n\le t\le 8n+3), \\ 
1\ge x(t)\ge -1 & (\mbox{if }\; 8n+3\le t\le 8n+4), \\ 
x(t)=-1 &  (\mbox{if }\; 8n+4\le t\le 8n+7), \\ 
-1\le x(t)\le 1 & (\mbox{if }\; 8n+7\le t\le 8n+8). 
\end{array}
\right.
\] 
%for any integer $n$.    
Let $y: \mathbb{R}\to \mathbb{R}$ be a $C^\infty$ 
periodic function of period $8$ satisfying 
the following condition for any integer $n$. 
\[
\left\{
\begin{array}{cc}
y(t)=-1 & (\mbox{if }\; 8n-2\le t\le 8n+1), \\ 
-1\le y(t)\le 1 & (\mbox{if }\; 8n+1\le t\le 8n+2), \\ 
y(t)=1 &  (\mbox{if }\; 8n+2\le t\le 8n+5), \\ 
1\ge x(t)\ge -1 & (\mbox{if }\; 8n+5\le t\le 8n+6).
\end{array}
\right.
\] 
%for any integer $n$.    
Define $\widetilde{f}: \mathbb{R}\to \mathbb{R}^2$ by 
$\widetilde{f}(t)=(x(t), y(t))$.    
Then, $\widetilde{f}$ is a $C^\infty$ periodic 
mapping of period $8$  
and the set of singular points of $\widetilde{f}$ contains 
infinitely many closed intervals  
\[
\cdots, [-2, -1], [0, 1], [2, 3], [4,5], [6,7], [8,9], \cdots.
\]  
{\color{black}From Figure \ref{figure 3}, it is clear that }
even if $\widetilde{f}$ has other singular points, 
the image of $\widetilde{f}$ is always the square with 
the following $4$ vertexes 
\begin{eqnarray*}
(1, -1)=\widetilde{f}([8n, 8n+1]),&  
(1,1)=\widetilde{f}([8n+2, 8n+3]), \\ 
(-1, 1)=\widetilde{f}([8n+4, 8n+5]),&  
(-1,-1)=\widetilde{f}([8n+6, 8n+7]).
\end{eqnarray*}
%%%%%%%%%%%%%%%%%%%%%%%%%%%%%%%%%%%%%%%%%%%%%%   
\begin{figure}
\begin{center}
\includegraphics[width=6cm, height=4cm]
{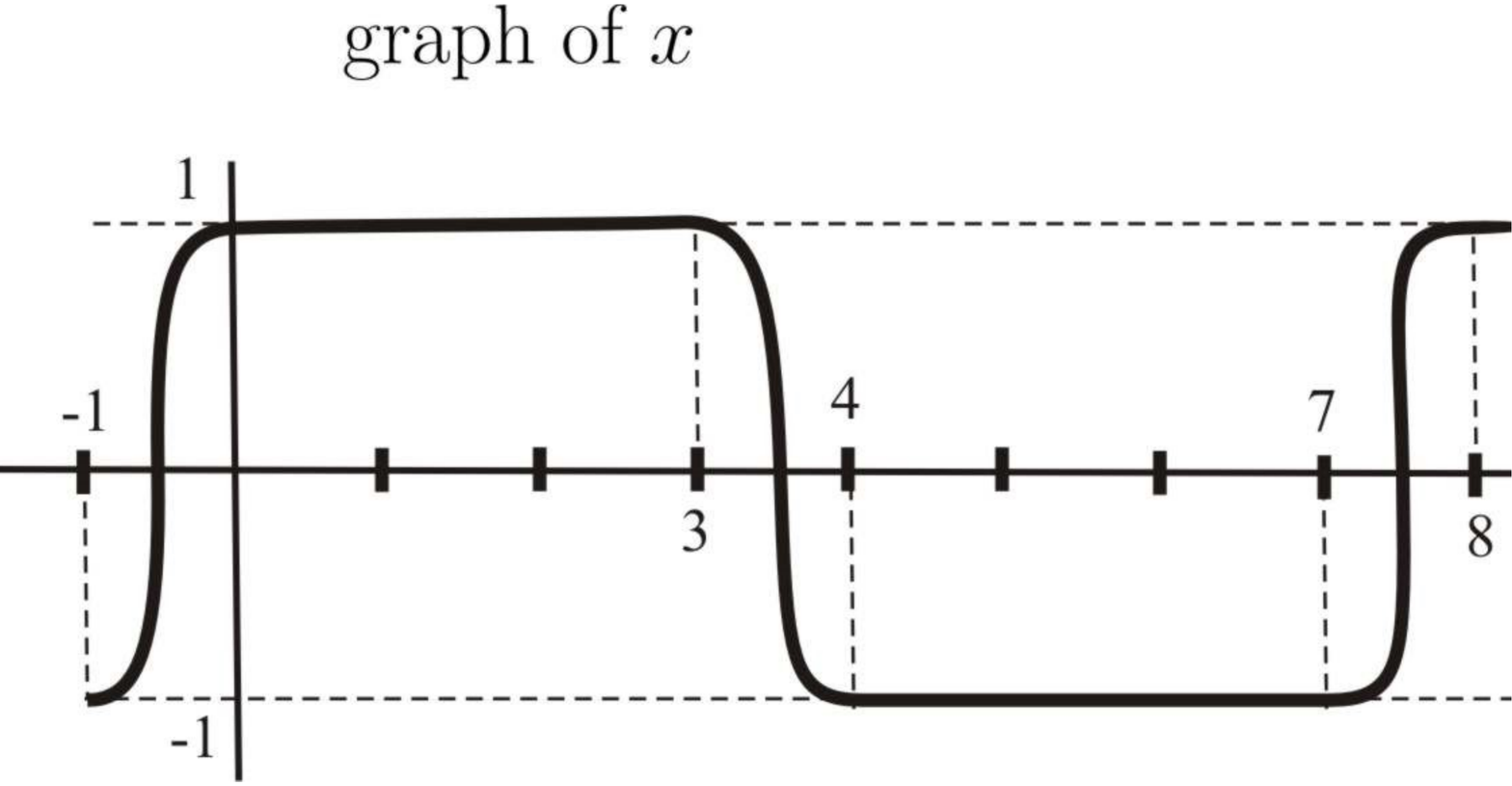}
\qquad 
\includegraphics[width=6.2cm, height=4cm]
{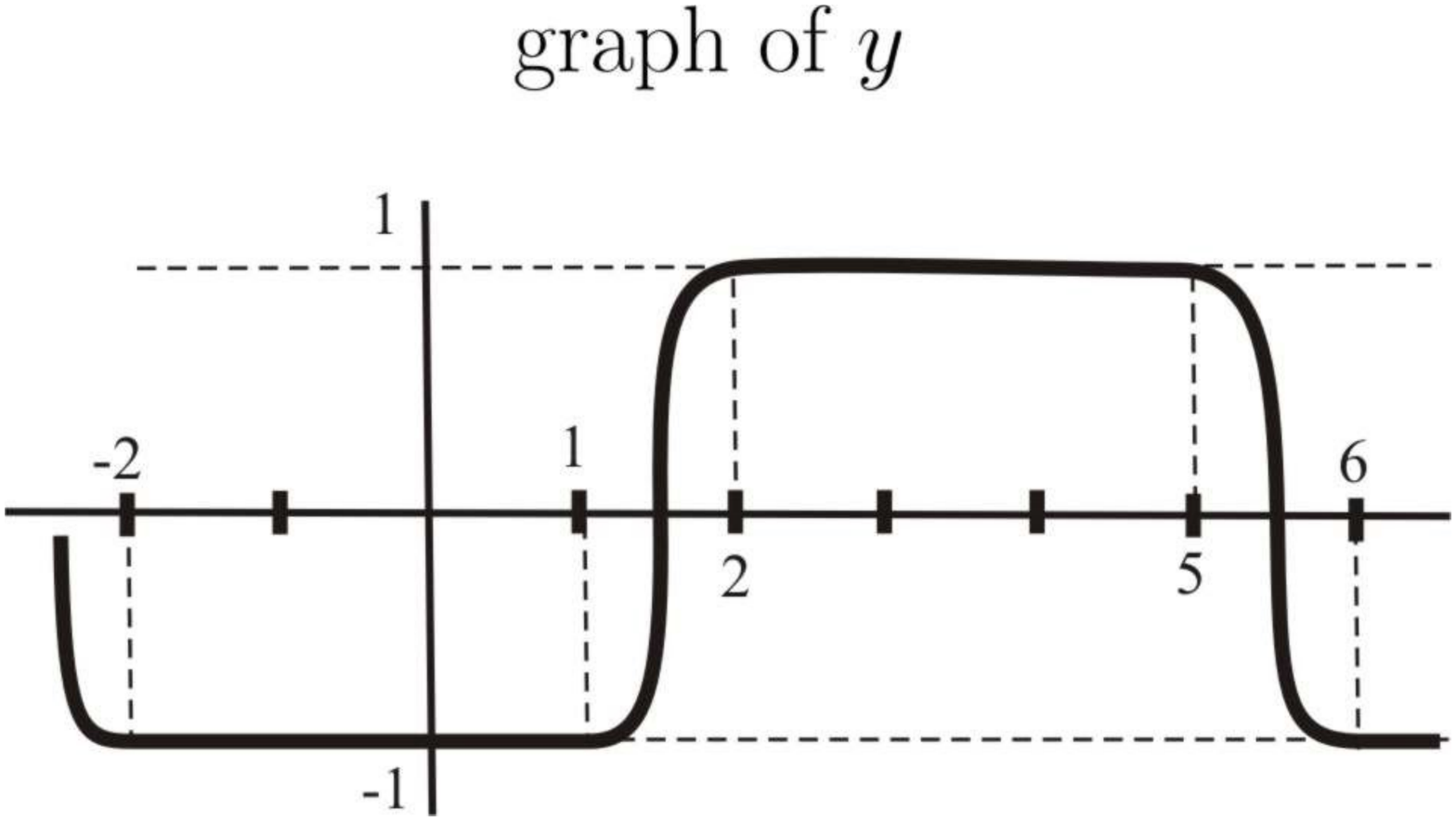}
\caption{%Left:
Graphs of $x$ and $y$. 
% Right: 
}
\label{figure 3}
\end{center}
\end{figure}
%%%%%%%%%%%%%%%%%%%%%%%%%%%%%%%%%%%%%%%%%%%%%%%
%%%%%%%%%%%%%%%%%%%%%%%%%%%%%%%%%%% 
%%%%%%%%%%%%%%%%%%%%%%%%%%%%%%%%%%% 
%\vspace{7cm}
%%%%%%%%%%%%%%%%%%%%%%%%%%%%%%%%%%% 
%%%%%%%%%%%%%%%%%%%%%%%%%%%%%%%%%%% 
\par 
Next, in order to assert that $\widetilde{f}$ is a frontal, 
we construct a non-zero normal vector 
$({\color{black}n_1}(t), {\color{black}n_2}(t))$ to $\widetilde{f}$ 
at $\widetilde{f}(t)$.     
Let ${\color{black}n_1}: \mathbb{R}\to \mathbb{R}$ be a $C^\infty$ 
periodic function of period $8$ satisfying 
the following condition for any integer $n$. 
\[
\left\{
\begin{array}{cc}
0\le {\color{black}n_1}(t)\le 1 
& (\mbox{if }\; 8n\le t\le 8n+1), \\ 
{\color{black}n_1}(t)=1 
& (\mbox{if }\; 8n+1\le t\le 8n+2), \\ 
1\ge {\color{black}n_1}(t)\ge 0 &  
(\mbox{if }\; 8n+2\le t\le 8n+3),  \\ 
{\color{black}n_1}(t)=0 
& (\mbox{if }\; 8n+3\le t\le 8n+4), \\  
0\ge {\color{black}n_1}(t)\ge -1 
& (\mbox{if }\; 8n+4\le t\le 8n+5), \\ 
{\color{black}n_1}(t)=-1 
& (\mbox{if }\; 8n+5\le t\le 8n+6), \\ 
-1\le {\color{black}n_1}(t)\le 0 
&  (\mbox{if }\; 8n+6\le t\le 8n+7), \\ 
{\color{black}n_1}(t)=0 
& (\mbox{if }\; 8n+7\le t\le 8n+8). \\ 
\end{array}
\right.
\] 
%for any integer $n$.                  
Let ${\color{black}n_2}: \mathbb{R}\to \mathbb{R}$ 
be a $C^\infty$ 
periodic function of period $8$ satisfying 
the following condition for any integer $n$. 
\[
\left\{
\begin{array}{cc}
-1\le {\color{black}n_2}(t)\le 0 
& (\mbox{if }\; 8n\le t\le 8n+1), \\ 
{\color{black}n_2}(t)=0 
& (\mbox{if }\; 8n+1\le t\le 8n+2), \\ 
0\le {\color{black}n_2}(t)\le 1 
&  (\mbox{if }\; 8n+2\le t\le 8n+3), \\ 
{\color{black}n_2}(t)=1 
& (\mbox{if }\; 8n+3\le t\le 8n+4), \\  
1\ge {\color{black}n_2}(t)\ge 0 
& (\mbox{if }\; 8n+4\le t\le 8n+5), \\ 
{\color{black}n_2}(t)=0 
& (\mbox{if }\; 8n+5\le t\le 8n+6), \\ 
0\ge {\color{black}n_2}(t)\ge -1 
&  (\mbox{if }\; 8n+6\le t\le 8n+7), \\ 
{\color{black}n_2}(t)=-1 
& (\mbox{if }\; 8n+7\le t\le 8n+8). \\ 
\end{array}
\right.
\] 
%for any integer $n$.           
{\color{black}From Figure \ref{figure 3} and 
Figure 4, %\ref{figure 4}, %\ref{figure4}, 
it is easily seen that %the following hold: 
the following two properties hold 
for any $t\in \mathbb{R}$.}    
\[ 
({\color{black}n_1}(t), {\color{black}n_2}(t))\ne (0,0) %\quad 
%(\forall t\in \mathbb{R})
\quad \mbox{ and }\quad  
\left(\frac{dx}{dt}(t), \frac{dy}{dt}(t)\right)
\cdot 
({\color{black}n_1}(t), {\color{black}n_2}(t))=0.      
\]   
Set 
\[
\widetilde{\nu}(t)=
\frac{1}{\sqrt{({\color{black}n_1}(t))^2+({\color{black}n_2}(t))^2}}
({\color{black}n_1}(t), {\color{black}n_2}(t)).   
\]
Then, $(\widetilde{f}, 
\widetilde{\nu}): \mathbb{R}\to \mathbb{R}^2\times S^1$ 
is actually a frontal.   
%%%%%%%%%%%%%%%%%%%%%%%%%%%%%%%%%%%%%%%%%%%%%%   
\begin{figure}
\begin{center}
\includegraphics[width=5.8cm, height=4cm]
{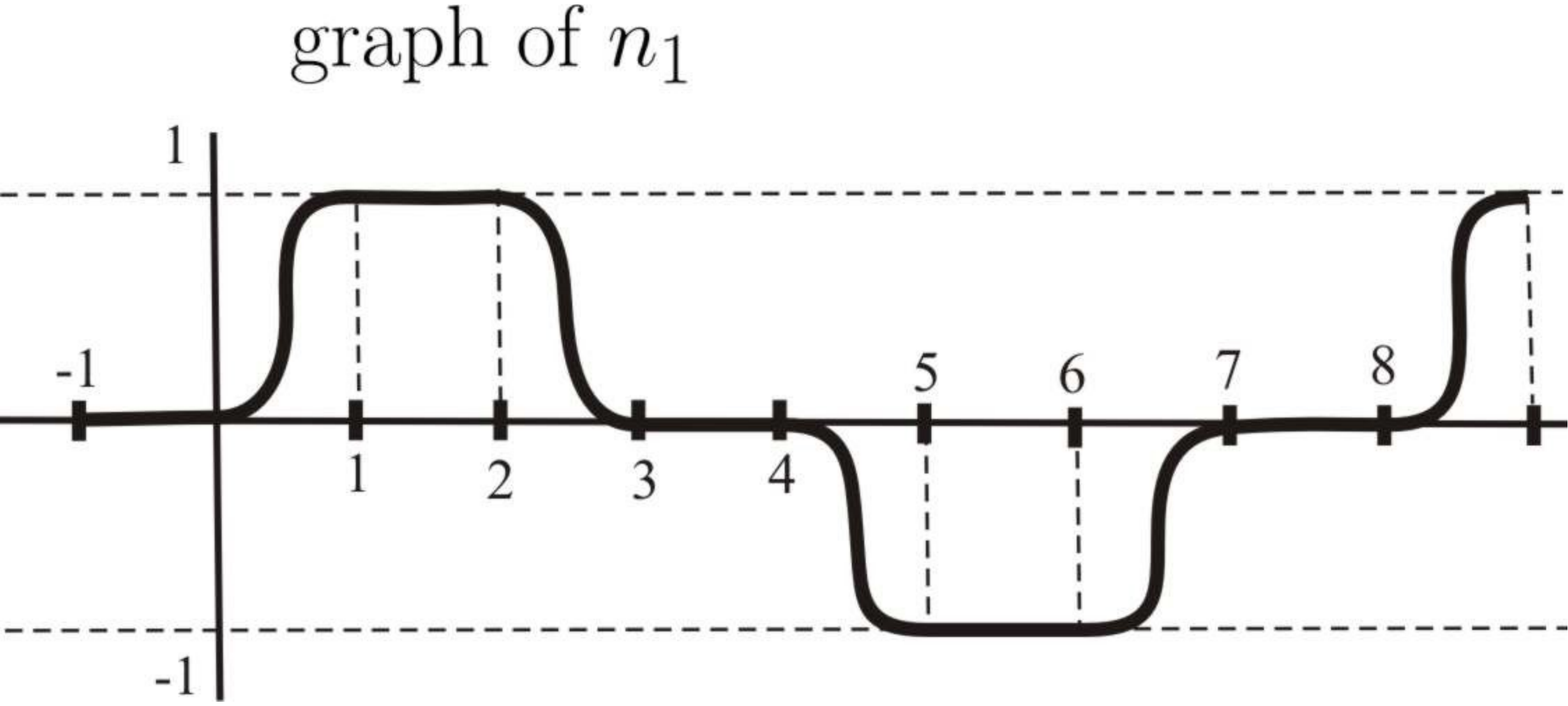}
\qquad 
\includegraphics[width=5.8cm, height=4cm]
{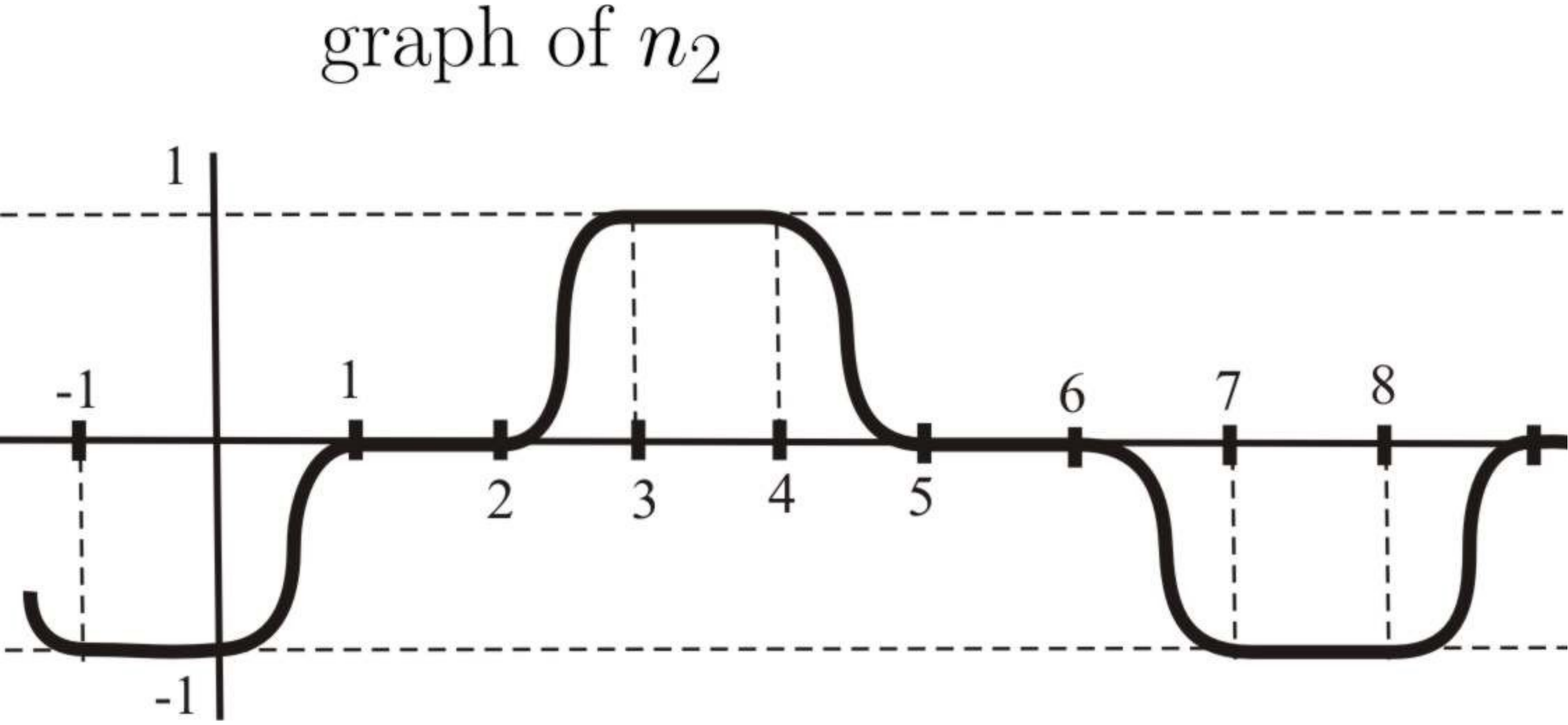}
{\color{black}\caption{%Left:
Graphs of $n_1$ and $n_2$.} 
% Right: 
}
\label{figure 4}
\end{center}
\end{figure}
%%%%%%%%%%%%%%%%%%%%%%%%%%%%%%%%%%%%%%%%%%%%%%%
%%%%%%%%%%%%%%%%%%%%%%%%%%%%%%%%%%%%%%%%% 
%%%%%%%%%%%%%%%%%%%%%%%%%%%%%%%%%%%%%%%%% 
%\vspace{7cm}
%%%%%%%%%%%%%%%%%%%%%%%%%%%%%%%%%%%%%%%%% 
%%%%%%%%%%%%%%%%%%%%%%%%%%%%%%%%%%%%%%%%% 
\par 
For the frontal $(\widetilde{f},\widetilde{\nu})$, 
the envelope of the one-parameter family 
of lines $\ell_t$ perpendicular to the unit vector 
$\widetilde{\nu}(t)$ and 
passing through the point $\widetilde{f}(t)$ 
does not restore the square 
$\widetilde{f}(\mathbb{R})$.    
\end{example}
These two examples show that, for frontals, the 
classical notion of envelope cannot restore 
the original figure in general.  
In order to eliminate the influence of singularities of frontals, 
Masatomo Takahashi has succeeded to 
improve the notion 
of envelopes (\cite{takahashi1, takahashi2}).      
The improvement due to Takahashi 
is nice, and thus for Example \ref{example 1}, 
the original figure $g(\mathbb{R})=S^1$ 
can be actually obtained 
as the envelope of Takahashi's sense.      
%Takahashi's envelope.    %for Example \ref{example 1}.   
However, unfortunately, the {\it variability condition} 
defined in \cite{takahashi1, takahashi2} is not satisfied for the 
frontal $(\widetilde{f}, \widetilde{\nu})$ given in 
Example \ref{example 2}.   
Thus, even Takahashi's envelope cannot restore 
the original square $\widetilde{f}(\mathbb{R})$ 
of Example \ref{example 2}.     
In Ishikawa's words, a frontal satisfying Takahashi's variability 
condition is called a 
{\it proper frontal} (\cite{ishikawa}, \S 6).     
The frontal 
$(\widetilde{f}, \widetilde{\nu}): 
\mathbb{R}\to \mathbb{R}^2\times S^1$ of  
Example \ref{example 2} is not a proper frontal.     
To the best of authors' knowledge, 
{\color{black}except for Example 2.5 given in \cite{ishikawa}}, 
all frontals investigated in detail so far  
are proper frontals.   
We would like to assert that non-proper frontals, too,  
are useful especially in application to surface science  
(see \S \ref{section 5}).          
\par     
The following Example \ref{example 3} shows 
that the uniqueness of anti-orthotomic (resp., negative pedal)    
does not hold in general even when the given mapping  
$f$ (resp., $g$) is a frontal.   
%although the orthotomic 
%of a given frontal relative to 
%a given point is always unique.  
\begin{example}\label{example 3} 
Let $f, g: \mathbb{R}\to \mathbb{R}^2$ be the constant 
mappings defined by $f(t)=(0,-1)$, $g(t)=(0,0)$.   
Let $\nu: \mathbb{R}\to S^1$ be the constant mapping 
defined by $\nu(t)=(-1, 0)$.    Set 
$P=(0,1)$ and define the 
constant mapping $\widetilde{\nu}: \mathbb{R}\to S^1$ 
by $\widetilde{\nu}(t)=(0,1)$.        
Then, for any $C^\infty$ mapping 
$\widetilde{f}: \mathbb{R}\to \mathbb{R}^2$ with the form 
$\widetilde{f}(t)=(\widetilde{f}_1(t), 0)$, 
the frontal $(\widetilde{f}, \widetilde{\nu})$ is  
an anti-orthotomic of $f$ relative to $P$ and 
a negative pedal of $g$ relative to $P$.   
\end{example}
%%%%%%%%%%%%%%%%%%%%%%%%%%%%%%%%%%%%%% 
%%%%%%%%%%%%%%%%%%%%%%%%%%%%%%%%%%%%%% 
The main purpose of this paper is 
to obtain the unique solution of 
the inverse problem for 
a given orthotomic $f$ relative to $P$ such that 
$P\in \mathcal{NS}_f$.   
%under reasonable conditions.  
%Namely, under reasonable conditions, 
%we obtain an explicit formula of an  
%anti-orthotomic for a given frontal 
%$f: N\to \mathbb{R}^{n+1}$ and a 
%point $P$.
\begin{theorem}\label{theorem 1}
Let $(f, \nu) : N\to \mathbb{R}^{n+1}\times S^n$, 
$P$ be a frontal 
%with 
%its unit normal ${\nu}: N\to \mathbb{R}^{n+1}$ 
and     
%satisfying conditions of frontal for ${f}$.   
a point of 
$\mathbb{R}^{n+1}$ such that 
$P\in \mathcal{NS}_f$ respectively.    
%Then, the orthotomic of $f$ relative to $P$ is a frontal.     
Let $\widetilde{f}: N\to \mathbb{R}^{n+1}$ be 
the mapping defined by 
\[
\widetilde{f}(x)=f(x)-\frac{||f(x)-P||^2}{2(f(x)-P)\cdot 
\nu(x)}\nu(x).    
\]
Then, the following four holds:   
\begin{enumerate}
\item[(1)] The mapping $\widetilde{f}$ 
is a frontal with its Gauss mapping 
$\widetilde{\nu}(x)=\frac{f(x)-P}{||f(x)-P||}$.   
\item[(2)] The mapping $\widetilde{f}$ is 
the unique anti-orthotomic of $f$ relative to $P$.   
\item[(3)] The point $P$ belongs to 
$\mathcal{NS}_{\widetilde{f}}$.   
\item[(4)] The equality 
$||\widetilde{f}(x)-P||=||\widetilde{f}(x)-f(x)||$ holds 
for any $x\in N$.   
%Moreover, the point $P$ belongs to the complement of 
%the image $\widetilde{f}(N)$.
\end{enumerate}
\end{theorem} 
%%%%%%%%%%%%%%%%%%%%%%%%%%%%%%%%% 
{\color{black}
In the case that $P=O$ and 
$f: I\to \mathbb{R}^2$ is a regular curve 
such that $f(x)\cdot \nu(x)\ne 0$ for any $x\in I$, 
the same formula for $\widetilde{f}$ has been already 
given in \cite{brucegiblin} 7.14 
as the envelope of perpendicular bysectors 
of segments  joining $f(x)$ and $P$.     
On the other hand, by Proposition 
\ref{proposition 1}, $\widetilde{f}(x)$ must be in the 
normal line $\{f(x)+a\nu(x)\; |\; a\in \mathbb{R}\}$.     
Therefore, in Theorem \ref{theorem 1}, 
just by solving simultaneous linear equations, 
the same formula for $\widetilde{f}$ can be obtained 
easily as the intersections of the perpendicular bysectors 
and the normal lines; and the no-silhouette condition 
$P\in \mathcal{NS}_f$ guarantees that 
each simultaneous linear equation must have 
the unique solution.   
} 
%%%%%%%%%%%%%%%%%%%%%%%%%%%%%%%%%
\par 
The following corollary clearly follows from 
Theorem \ref{theorem 1}  
\begin{corollary}\label{corollary 2}
Let $(g, \nu) : N\to \mathbb{R}^{n+1}\times S^n$, 
$P$ be a frontal 
%with 
%its unit normal ${\nu}: N\to \mathbb{R}^{n+1}$ 
and     
%satisfying conditions of frontal for ${f}$.   
a point of 
$\mathbb{R}^{n+1}$ such that 
$P\in \mathcal{NS}_g$ respectively.    
%Then, the orthotomic of $f$ relative to $P$ is a frontal.     
Let $\widetilde{f}: N\to \mathbb{R}^{n+1}$ be 
the mapping defined by 
\[
\widetilde{f}(x)=2g(x)-P-
\frac{||g(x)-P||^2}{(g(x)-P)\cdot 
\nu(x)}\nu(x).    
\]
Then, the following four holds:   
\begin{enumerate}
\item[(1)] The mapping $\widetilde{f}$ 
is a frontal with its Gauss mapping  
$\widetilde{\nu}(x)=\frac{g(x)-P}{||g(x)-P||}$.   
\item[(2)] The mapping $\widetilde{f}$ is 
the unique negative pedal of $g$ relative to $P$.   
\item[(3)] The point $P$ belongs to 
$\mathcal{NS}_{\widetilde{f}}$.   
\item[(4)] The equality 
$||\widetilde{f}(x)-P||=||\widetilde{f}(x)-2g(x)+P||$ holds 
for any $x\in N$.   
%Moreover, the point $P$ belongs to the complement of 
%the image $\widetilde{f}(N)$.
\end{enumerate}
\end{corollary} 
%Any embedding 
%$\widetilde{f}: \mathbb{R}^n\to \mathbb{R}^n\times \{0\}$ 
%given in Example \ref{example 1} may be constructed 
%by Theorem \ref{theorem 1}.      
%Moreover, if $f: \mathbb{R}^n\to \mathbb{R}^{n+1}$ 
%(resp., $f: S^n\to \mathbb{R}^{n+1}$) 
%is a parametrization of a hyperplane (resp., the inclusion), 
%then for any point $P$ of the complement of its image, 
%$\widetilde{f}$ defined in Theorem \ref{theorem 1} is 
%nothing but 
%a parametrization of the parabolic hypersurface 
%with the focus $P$ and the directrix $f$ 
%(resp., the elliptic hypersurface with focuses $P$ and the %origin $O$ 
%defined by $||\widetilde{f}(x)-P||+||\widetilde{f}(x)-O||=1$).     
%For details, see Example \ref{hyperplane} in Section %\ref{section 5}.   
%%%%%%%%%%%%%%%%%%%%%%%%%%%%%%%%%%%%%%%%%%%%%%%%%%%%%%%%%%%%   
%%%%%%%%%%%%%%%%%%%%%%%%%%%%%%%%%%%%%%%%%%%%%%%%%%%%%%%%%%%% 
\par 
\bigskip 
This paper is organized as follows. 
In Section \ref{section 2}, 
Proposition \ref{proposition 1} is proved.   
Theorem \ref{theorem 1} is proved in Section 
\ref{section 3}.    
%By definition, the results  
%on pedals and negative pedals corresponding to  
%Proposition \ref{proposition 1} and Theorem \ref{theorem 1} 
%can be obtained immediately.   
%Section \ref{section 4} devotes to collect such results.     
In the case that the Gauss 
mapping $\widetilde{\nu}$ of $\widetilde{f}$ is the identity 
mapping,  
there is the famous Cahn-Hoffman vector formula for 
$\widetilde{f}$ (\cite{hoffmancahn}).   
%This formula is proved without using envelopes.    
In Section \ref{section 4}, as the first application of 
Theorem \ref{theorem 1}, Cahn-Hoffman formula  
is shown. 
% by Theorem \ref{theorem 1}.    
In Section \ref{section 5}, as the second application 
of Theorem \ref{theorem 1}, the optical meaning of the 
anti-orthotomic $\widetilde{f}$ is clarified 
even at a singular point of the Gauss mapping 
$\widetilde{\nu}$ of $\widetilde{f}$.   
Moreover, in order to show how 
the clarified optical meaning 
is useful, it is applied to construct the exact shape of 
the orthotomic $f$ for the frontal $\widetilde{f}$ in Example 
\ref{example 2} and a given point 
$P\in \mathcal{NS}_{\widetilde{f}}$.    
Finally, in Section \ref{section 6}, as the third application 
of Theorem \ref{theorem 1}, 
it is given 
a criterion 
that a given frontal is actually a front.      
%
%Theorem \ref{theorem 1} gives an exact solution for 
%an inverse problem.   By Proposition 
%\ref{proposition 1} and Theorem \ref{theorem 1}, 
%one can manipulate normal vectors of two kinds as one like.      
%Thus, we believe that 
%Theorem \ref{theorem 1} must have many applications.    
%In \S\S 4--6, applications of Theorem \ref{theorem 1} 
%obtained so far are given.   
%   
%Finally, in Section \ref{section 5}, 
%iterated reflections are investigated.    
%%%%%%%%%%%%%%%%%%%%%%%%%%%%%%%%%
%%%%%%%%%%%%%%%%%%%%%%%%%%%%% 
%%%%%%%%%%%%%%%%%%%%%%%%%%%%%%%%
%%%%%%%%%%%%%%%%%%%%%%%%%%%%%% 
\section{Proof of Proposition \ref{proposition 1}}
\label{section 2}
%%%%%%%%%%%%%%%%%%%%%%%%%%%%%%%%%%%%%%%%
%%%%%%%%%%%%%%%%%%%%%%%%%%%%%%%%%%%%%%%%
\subsection{Proof that $f$ is a frontal 
with its Gauss mapping 
$\nu(x)=
\frac{f(x)-\widetilde{f}(x)}{||f(x)-\widetilde{f}(x)||}$}
%%%%%%%%%%%%%%%%%%%%%%%%%%%%%%%%%%%%%%%%%%%%%%%%% 
\quad \\ 
Recall that $f$ is defined by 
$f(x)=2\left((\widetilde{f}(x)-P)\cdot 
\widetilde{\nu}(x)\right)\widetilde{\nu}(x)+P$.    
\begin{lemma}\label{lemma 1}
For any $x\in N$, $f(x)-\widetilde{f}(x)$ is a non-zero vector.  
\end{lemma}
\proof 
Suppose that $f(x_0)=\widetilde{f}(x_0)$ for some $x_0\in N$.    
Then, for the $x_0$, the following holds: 
\[
\widetilde{f}(x_0)-P=2\left((\widetilde{f}(x_0)-P)\cdot 
\widetilde{\nu}(x_0)\right)\widetilde{\nu}(x_0).   
\]
This implies $
%\begin{enumerate}
%\item[(1)] There exists a real number $a$ such that 
%$f(x_0)-P=a\nu(x_0)$.   
(\widetilde{f}(x_0)-P)\cdot \widetilde{\nu}(x_0)=0,    
$ 
which contradicts the assumption  
$P\in \mathcal{NS}_{\widetilde{f}}$. 
\hfill 
$\Box$ 
\par 
\medskip    
Set 
\[
\nu(x)=\frac{f(x)-\widetilde{f}(x)}{||f(x)-\widetilde{f}(x)||}.
\]   
Then, 
it is sufficient to show that 
$d{f}_x({\bf v})\cdot {\nu}(x)=0$ for any 
$x\in N$ and any ${\bf v}\in T_x N$.   
In other words, it is sufficient to show that 
\[
(f\circ \xi)'(0)\cdot \nu(x)=0   
\]
for any curve 
$\xi: (-\varepsilon, \varepsilon)\to N$ such that $\xi(0)=x$.     
The following lemma clearly holds:   
\begin{lemma}\label{lemma 2}    
\begin{enumerate}
\item[(1)] $(\widetilde{f}\circ \xi)'(0)\cdot 
\widetilde{\nu}(x)=0$.  
\item[(2)] $(\widetilde{\nu}\circ \xi)'(0)\cdot 
\widetilde{\nu}(x)=0$.   
\item[(3)] 
$\left(f\circ \xi\right)'(0) = 
2\left(\left(\widetilde{f}(x)-P\right)\cdot 
\left(\widetilde{\nu}\circ\xi\right)'(0)\right)
\widetilde{\nu}(x)+ 
2\left(\left(\widetilde{f}(x)-P\right)\cdot 
\widetilde{\nu}(x)\right)
\left(\widetilde{\nu}\circ \xi\right)'(0).   
$
\end{enumerate}
\end{lemma}
By using Lemma \ref{lemma 2}, we have the following:   
\begin{eqnarray*}
{ } & { } & ||f(x)-\widetilde{f}(x)||
\left(\left(f\circ \xi\right)'(0)\cdot 
\nu(x)\right) \\ 
{ } & = & \left(f\circ \xi\right)'(0)\cdot 
\left(f(x)-\widetilde{f}(x)\right) \\ 
{ } & = & 
\left(f\circ \xi\right)'(0)\cdot 
\left(2\left(\left(\widetilde{f}(x)-P\right)\cdot 
\widetilde{\nu}(x)\right)\widetilde{\nu}(x)
-\left(\widetilde{f}(x)-P\right)\right) \\ 
{ } & = & 
4\left((\widetilde{f}(x)-P)\cdot 
\left(\widetilde{\nu}\circ \xi\right)'(0)\right)
\left((\widetilde{f}(x)-P)\cdot \widetilde{\nu}(x)\right) \\ 
{ } & { } & 
\qquad 
-2\left((\widetilde{f}(x)-P)\cdot 
\left(\widetilde{\nu}\circ \xi\right)'(0)\right) 
\left(\widetilde{\nu}(x)\cdot (\widetilde{f}(x)-P)\right)
\\ 
{ } & { } & 
\qquad\qquad  
-2\left((\widetilde{f}(x)-P)\cdot \widetilde{\nu}(x)\right) 
\left(\left(\widetilde{\nu}\circ \xi\right)'(0)\cdot 
(\widetilde{f}(x)-P)\right)
\\ 
{ } & = & 0.      
\end{eqnarray*}
\hfill 
$\Box$ 
%\begin{note}\label{note 1}
%For the proof given in this subsection, the assumption 
%$P\in \mathcal{NS}_f$ 
%is too strong.   
%The weaker condition $P\not\in f(N)$ is enough.     
%However, for the proof given in the next subsection, 
%the assumption must be 
%$P\in \mathcal{NS}_f$.       
%\end{note}
%%%%%%%%%%%%%%%%%%%%%
%%%%%%%%%%%%%%%%%%%%%%%%%%%%%%% 
\subsection{Proof that 
$\left(f(x)-P\right)\cdot \nu(x)\ne 0$ 
for any $x\in N$}    
%of $P\not\in \widetilde{f}(N)$}
%Recall the following:  
%\begin{eqnarray*}
%\widetilde{f} & = & 2\left((f(x)-P)\cdot \nu(x)\right)\nu(x)+P 
%\\ 
%\widetilde{\nu} & = & \frac{\widetilde{f}(x)-f(x)}
%{||\widetilde{f}(x)-f(x)||}.   
%\end{eqnarray*}
\quad 
\\
For any $x\in N$, we have the following:   
\begin{eqnarray*}
{ } & { } & 
||f(x)-\widetilde{f}(x)||
\left(f(x)-P\right)\cdot \nu(x) \\ 
{ } & = & 
\left(f(x)-P\right)\cdot 
\left(f(x)-\widetilde{f}(x)\right) \\ 
{ } & = & 
2\left(\left(\left(\widetilde{f}(x)-P\right)\cdot 
\widetilde{\nu}(x)\right)\widetilde{\nu}(x)\right)\cdot 
\left(f(x)-\widetilde{f}(x)\right) \\ 
{ } & = & 
2\left(\left(\left(\widetilde{f}(x)-P\right)\cdot 
\widetilde{\nu}(x)\right)\widetilde{\nu}(x)\right)\cdot 
\left(2\left(\left(\widetilde{f}(x)-P\right)
\cdot \widetilde{\nu}(x)\right)
\widetilde{\nu}(x)- \left(\widetilde{f}(x)-P\right)\right) \\ 
{ } & = & 
4\left(\left(\widetilde{f}(x)-P\right)
\cdot \widetilde{\nu}(x)\right)^2
-2\left(\left(\widetilde{f}(x)-P\right)\cdot 
\widetilde{\nu}(x)\right) 
\left(\widetilde{\nu}(x)\cdot 
\left(\widetilde{f}(x)-P\right)\right)
\\ 
{ } & = & 
2\left(\left(\widetilde{f}(x)-P\right)\cdot 
\widetilde{\nu}(x)\right)^2.
\end{eqnarray*}
By the assumption $P\in \mathcal{NS}_{\widetilde{f}}$, 
it follows that $\left(f(x)-P\right)\cdot \nu(x)
\ne 0$ for any $x\in N$.    
\hfill 
$\Box$
%%%%%%%%%%%%%%%%%%%%%%%%%%%%%%%%%%%%%%%%%%%%%%%%% 
%%%%%%%%%%%%%%%%%%%%%%%%%%%%%%%%%%%%%%%%%%%%%%%%%
%%%%%%%%%%%%%%%%%%%%%%%%%%%%%%%%%%%%%%%%%%%%%%%%%
\section{Proof of Theorem \ref{theorem 1}}\label{section 3}
%%%%%%%%%%%%%%%%%%%%%%%%%%%%%%%%%%%%%%%%%%%%%%%%% 
\subsection{Proof that $\widetilde{f}$ is a frontal 
with Gauss mapping 
$\widetilde{\nu}(x)=\frac{f(x)-P}{||f(x)-P||}$}
%%%%%%%%%%%%%%%%%%%%%%%%%%%%%%%%%%%%%%%%%%%%%%%%% 
\quad \\ 
From the assumption that $(f(x)-P)\cdot \nu(x)\ne 0$ for any 
$x\in N$, it follows that $f(x)\ne P$ for any $x\in N$.   
Thus,  
\[
\widetilde{\nu}(x)=\frac{f(x)-P}{||f(x)-P||}.    
\]
is well-defined. 
%a unit vector for any $x\in N$.
Then, it is sufficient to show that 
\[
(\widetilde{f}\circ \xi)'(0)\cdot \widetilde{\nu}(x)=0   
\]
for any curve $\xi: (-\varepsilon, \varepsilon)\to N$ 
such that $\xi(0)=x$.  
Since 
$\widetilde{f}$ has the form 
\[
\widetilde{f}(x)=f(x)-\frac{||f(x)-P||^2}{2\left(f(x)-P\right)
\cdot \nu(x)}\nu(x), 
\] 
we have the following:   
\begin{eqnarray*}
{ } & { } & 
||f(x)-P||
\left(\left(\widetilde{f}\circ \xi\right)'(0)\cdot 
\widetilde{\nu}(x)\right) \\ 
{ } & = &  
\left(f\circ \xi\right)'(0)\cdot \left(f(x)-P\right) 
-\frac{\left(f\circ \xi\right)'(0)\cdot \left(f(x)-P\right)}
{\left(f(x)-P\right)\cdot \nu(x)}
\left(\nu(x)\cdot \left(f(x)-P\right)\right) \\ 
{ } & { } & 
\qquad 
+\frac{\left(f(x)-P\right)\cdot\left(\nu\circ \xi\right)'(0)}
{2\left(\left(f(x)-P\right)\cdot \nu(x)\right)^2}||f(x)-P||^2
\left(\nu(x)\cdot \left(f(x)-P\right)\right) \\
{ } & { } & 
\qquad \qquad 
-\frac{||f(x)-P||^2}
{2\left(\left(f(x)-P\right)\cdot \nu(x)\right)}
\left(\left(\nu\circ \xi\right)'(0)\cdot 
\left(f(x)-P\right)\right) \\ 
{ } & = & 
\left(f\circ \xi\right)'(0)\cdot \left(f(x)-P\right) 
-\left(f\circ \xi\right)'(0)\cdot \left(f(x)-P\right) \\ 
{ } & { } & 
\qquad 
+\frac{\left(f(x)-P\right)\cdot\left(\nu\circ \xi\right)'(0)}
{2\left(\left(f(x)-P\right)\cdot \nu(x)\right)}||f(x)-P||^2  \\ 
{ } & { } & 
\qquad \qquad 
-\frac{||f(x)-P||^2}
{2\left(\left(f(x)-P\right)\cdot \nu(x)\right)}
\left(\left(\nu\circ \xi\right)'(0)\cdot 
\left(f(x)-P\right)\right) \\ 
{ } & = & 
0+0=0.   
\end{eqnarray*}
\hfill 
$\Box$
%%%%%%%%%%%%%%%%%%%%%%%%%%%%%%%%%%%%%%%%%%%%%%%%% 
\subsection{Proof that $\widetilde{f}$ is 
{\color{black}the unique} anti-orthotomic of $f$ 
relative to $P$}\label{subsection 3.2}
%%%%%%%%%%%%%%%%%%%%%%%%%%%%%%%%%%%%%%%%%%%%%%%%% 
\quad \\ 
%%%%%%%%%%%%%%%%%%%%%%%%%%%%%%%  
{\color{black}The proof is essentially 
given in the paragraph next to Theorem \ref{theorem 1}.   
Thus, in this subsection, just a confirmation by definition 
is given.   
}    
%%%%%%%%%%%%%%%%%%%%%%%%%%%% 
Recall that %in Subsection \ref{subsection 3.1} the mapping 
$\widetilde{\nu}(x)=\frac{f(x)-P}{||f(x)-P||}$.  
%: N\to \mathbb{R}^{n+1}$ was defined as follows 
We have the following:       
\begin{eqnarray*}
{ } & { } & 
2\left(\left(\widetilde{f}(x)-P\right)\cdot 
\widetilde{\nu}(x)\right)
\widetilde{\nu}(x)+P   \\ 
{ } & = & 
2\left(\left(
\left(f(x)-P\right)- \frac{||f(x)-P||^2}{2\left(f(x)-P\right)
\cdot \nu(x)}
\nu(x)
\right)\cdot 
\frac{f(x)-P}{||f(x)-P||}
\right)
\frac{f(x)-P}{||f(x)-P||} + P \\ 
{ } & = & 
2\left(||f(x)-P||-\frac{||f(x)-P||}{2}\right)
\frac{f(x)-P}{||f(x)-P||} + P \\ 
{ } & = & 
f(x)-P+P=f(x).
\end{eqnarray*}
\hfill 
$\Box$  
%%%%%%%%%%%%%%%%%%%%%%%%%%%%%%%%%%%%%%%%%%%%%%%%%%%% 
\subsection{Proof that 
$\left(\widetilde{f}(x)-P\right)\cdot \widetilde{\nu}(x)\ne 0$ 
for any $x\in N$}    \label{subsection 3.3}
%of $P\not\in \widetilde{f}(N)$}
%%%%%%%%%%%%%%%%%%%%%%%%%%%%%%%%%%%%%%%%%%%%%%%%% 
\quad \\ 
%Recall that $\widetilde{\nu}(x)=\frac{f(x)-P}{||f(x)-P||}$.   
For any $x\in N$ we have the following:  
\begin{eqnarray*}
\left(\widetilde{f}(x)-P\right)\cdot \widetilde{\nu}(x) 
& = & \left(\widetilde{f}(x)-P\right)\cdot 
\frac{f(x)-P}{||f(x)-P||} \\ 
{ } & = & 
\left(f(x)-\frac{||f(x)-P||^2}
{2(f(x)-P)\cdot \nu(x)}\nu(x)-P\right)
\cdot \frac{f(x)-P}{||f(x)-P||} \\ 
{ } & = & 
\left(\left(f(x)-P\right)-\frac{||f(x)-P||^2}
{2(f(x)-P)\cdot \nu(x)}\nu(x)\right)
\cdot \frac{f(x)-P}{||f(x)-P||} \\ 
{ } & = & 
||f(x)-P|| - \frac{||f(x)-P||}{2} \\ 
{ } & = & 
\frac{||f(x)-P||}{2}\ne 0.  
\end{eqnarray*}
%Suppose that $P=\widetilde{f}(x_0)$ for some $x_0\in N$.     
%Then, the following holds:   
%\[
%f(x_0)-P=\frac{||f(x_0)-P||^2}{2\left(f(x_0)-P\right)\cdot \nu(x_0)}\nu(x_0).  
%\leqno{(*)}
%\]
%Let $\theta_0$ be the angle between two non-zero vectors 
%$f(x_0)-P$ and $\nu(x_0)$.   Then, 
%the above equality $(*)$ implies that $\cos\theta_0$ is $1$ or 
%$-1$.    
%On the other hand, dotting $\nu(x_0)$ with both sides of $(*)$ 
%yields the following:    
%\[
%2\left(\left(f(x)-P\right)\cdot \nu(x_0)\right)^2=||f(x)-P||^2.   
%\]
%Thus, we have the following 
%which contradicts $\cos\theta_0=1$ or 
%$-1$.     
%\[
%2\cos^2\theta_0=1.    
%\]
%Therefore, $P$ must be outside $\widetilde{f}(N)$.   
\hfill 
$\Box$ 
%%%%%%%%%%%%%%%%%%%%%%%%%%%%%%%%%%%%%%%%%%%%%%%%%%% 
\subsection{Proof that the equality 
$||\widetilde{f}(x)-P||=||\widetilde{f}(x)-f(x)||$ holds 
for any $x\in N$
}
\label{subsection 3.4}
%%%%%%%%%%%%%%%%%%%%%%%%%%%%%%%%%%%%%%%%%%%%%%%%% 
Since $\widetilde{f}(x)-P
=(f(x)-P)-\frac{||f(x)-P||^2}{2(f(x)-P)\cdot \nu(x)}\nu(x)$, 
the following holds for any $x\in N$:   
\begin{eqnarray*}
{ } & { } & \left\|\widetilde{f}(x)-P\right\|^2 \\ 
{ } & = & 
\left(\widetilde{f}(x)-P\right)\cdot 
\left(\widetilde{f}(x)-P\right) \\ 
{ } & = & 
\left((f(x)-P)-\frac{||f(x)-P||^2}{2(f(x)-P)\cdot 
\nu(x)}\nu(x)\right)\cdot 
\left((f(x)-P)-\frac{||f(x)-P||^2}{2(f(x)-P)\cdot 
\nu(x)}\nu(x)\right) \\ 
{ } & = & 
||f(x)-P||^2 -||f(x)-P||^2+\frac{||f(x)-P||^2}{4\left((f(x)-P)
\cdot \nu(x)\right)^2} \\ 
{ } & = & 
\frac{||f(x)-P||^2}{4\left((f(x)-P)\cdot \nu(x)\right)^2} \\ 
{ } & = & 
\left\|\widetilde{f}(x)- f(x)\right\|^2.   
\end{eqnarray*}
\hfill 
$\Box$ 
%%%%%%%%%%%%%%%%%%%%%%%%%%%%%%%%%%%%%%%%%%%%%%%%% 
%%%%%%%%%%%%%%%%%%%%%%%%%%%%%%%%%%%%%%%%%%%%%%%%%%%% 
%%%%%%%%%%%%%%%%%%%%%%%%%%%%%%%%%%%%%%%%%%%%%%%%% 
\section{Application 1: 
Generalization of Cahn-Hoffman vector formula}
\label{section 4} 
%%%%%%%%%%%%%%%%%%%%%%%%%%%%%%%%%%%%%%%%%%%%%%%%% 
Let $(g, \nu): N\to \mathbb{R}^{n+1}\times S^n$ 
be a frontal.     We assume that $\mathcal{NS}_g$ 
is not empty.     Let $P$ be a point of $\mathcal{NS}_g$.   
Then, by Corollary \ref{corollary 2}, the mapping  
$\left(\widetilde{f}, \widetilde{\nu}\right): 
N\to \mathbb{R}^{n+1}\times S^n$ defined by 
\begin{eqnarray*}
\widetilde{f}(x) & = & 2g(x) - P - 
\frac{||g(x)-P||^2}{(g(x)-P)\cdot \nu(x)}\nu(x), \\ 
\widetilde{\nu}(x) & = & 
\frac{g(x)-P}{||g(x)-P||}
\end{eqnarray*}
is a frontal and the unique negative pedal of $g$ 
relative to $P$.    
Set $\gamma(x)=||g(x)-P||$.    
Then, by using $\widetilde{\nu}: N\to S^n$ 
and $\gamma: N\to \mathbb{R}_+$, 
$g(x)-P$ can be expressed as follows.      
\[
g(x)-P=\gamma(x)\widetilde{\nu}(x).    
\]
In \cite{hoffmancahn}, under the assumption that 
$N$ is the unit sphere $S^n$ 
and $\widetilde{\nu}: S^n\to S^n$ is the identity 
mapping and 
under the identification   
$\mathbb{R}^{n+1}=
T_{\widetilde{\nu}(x)}\mathbb{R}^{n+1}
=T_{\widetilde{\nu}(x)}S^n
\oplus \mathbb{R}\widetilde{\nu}(x)$,  
D.~W.~Hoffman and 
J.~W.~Cahn showed the following.   
\begin{theorem}
[Cahn-Hoffman vector formula  \cite{hoffmancahn}]
\label{hoffmancahn}
For any $x\in S^n$, the following equality holds.  
\[
\widetilde{f}(x)-g(x)=\nabla\gamma(x)
\oplus 0\widetilde{\nu}(x).   
\]
\end{theorem}
\noindent 
Here, $\nabla\gamma(x)$ stands for the gradient 
vector of $\gamma$ at $x$ with respect to 
the normal coordinate system of $S^n$ around $x$.    
In this section, as an application of 
Theorem \ref{theorem 1}, we generalize Theorem 
\ref{hoffmancahn} as follows.   
\begin{theorem}\label{theorem 2}
Let $(g, \nu): N\to \mathbb{R}^{n+1}\times S^n$ be a frontal 
and let $P$ be a point of $\mathcal{NS}_g$.        
%Let $x_0\in N$ and $P\in \mathbb{R}^{n+1}$ satisfy  
%$(g(x_0)-P)\cdot \nu(x_0)\ne 0$.   
Suppose that $\widetilde{\nu}$ is non-singular 
at $x$.   
Then,  the following equality holds.        
\[
\widetilde{f}(x) - g(x) 
=   
{\left(\left(J\widetilde{\nu}(x)\right)^{-1}\right)^t
\nabla \gamma(x)
{\oplus {0}\widetilde{\nu}(x)}}.   
\]
\end{theorem}
\noindent 
Here, $\left(\left(J\widetilde{\nu}(x)\right)^{-1}\right)^t$ 
stands for the transposed matrix of 
the inverse of Jacobian matrix of $\widetilde{\nu}$ 
with respect to an arbitrary local coordinate system 
around $x\in N$ and the normal coordinate system 
around $\widetilde{\nu}(x)\in S^n$.      
Theorem \ref{theorem 2} yields not only 
Theorem \ref{hoffmancahn} but also 
the following.   
\begin{corollary}   
Let $(g, \nu): N\to \mathbb{R}^{n+1}\times S^n$ be a frontal 
and let $P$ be a point of $\mathcal{NS}_g$.        
%Let $x_0\in N$ and $P\in \mathbb{R}^{n+1}$ satisfy  
%$(g(x_0)-P)\cdot \nu(x_0)\ne 0$.   
Suppose that $\widetilde{\nu}$ is non-singular 
at $x$.  
Then, $x$ is a singular point of $\gamma$ if and only 
if $\widetilde{f}(x)=g(x)$ is satisfied.    
%where $\gamma: N\to \mathbb{R}$ is a 
%$C^\infty$ function defined by $\gamma(x)=||g(x)-P||$.   
\end{corollary}
\noindent 
{\bf Proof of Theorem \ref{theorem 2}.} 
\par 
Since  
\[
\widetilde{f}(x)=f(x) 
- \frac{||f(x)-P ||^2}{2(f(x)-P)\cdot \nu(x)}\nu(x), 
\]  
it follows  
%of THEOREM 1, 
\begin{eqnarray*}
 f(x)-\widetilde{f}(x) 
 & = & 
\frac{||f(x)-P||^2}{2(f(x)-P)\cdot \nu(x)}\nu(x) \\ 
{} & = & 
\frac{4\gamma^2(x)}
{4\gamma(x)\left(\widetilde{\nu}(x)\cdot \nu(x)\right)}
\nu(x) \\
{} & = & 
\frac{\gamma(x)}
{\widetilde{\nu}(x)\cdot \nu(x)}
\nu(x) \\ 
{} & = & 
\frac{\gamma(x)}
{\widetilde{\nu}(x)\cdot \nu_2(x)}
\left(\nu_1(x)\oplus \nu_2(x)\right), 
\end{eqnarray*}
where $\nu(x)=\nu_1(x)\oplus \nu_2(x)$  
and 
$\nu_1(x)\in T_{\widetilde{\nu}(x)}S^n, 
\nu_2(x)\in \mathbb{R}\widetilde{\nu}(x)$.   
\par 
In order to represent $\nu_1, \nu_2$ in terms of 
$\gamma$ and $\widetilde{\nu}$, the same technique 
as in \cite{nishimura} is used.   
%On the other hand, 
Since 
$f(x)-P=2\gamma(x)\widetilde{\nu}(x)$, 
for any $v\in T_x N$, 
\[
df_x(v)=2\gamma(x)d\widetilde{\nu}_x(v)\oplus 
2d\gamma_x(v)\widetilde{\nu}(x), 
\]
where 
$\gamma(x)d\widetilde{\nu}_x(v)\in T_{\widetilde{\nu}(x)}S^n$ 
and $d\gamma_x(v)\widetilde{\nu}(x)\in 
\mathbb{R}\widetilde{\nu}(x)$.  
Thus, the Jacobian matrix $Jf$ of $f$ at $x$ 
with respect to an arbitrary local coordinate system 
around $x\in N$ and the direct product 
of the normal coordinate system and $\mathbb{R}$ 
around $f(x)=2\gamma(x)\widetilde{\nu}(x)$ has the 
following form, where $J\widetilde{\nu}(x)$ stands 
for the Jacobian matrix of $\widetilde{\nu}$ at $x$ 
and $\left(\nabla\gamma(t)\right)^t$ stands for the 
transposed vector of the gradient of $\gamma$ at $x$.   
\[
Jf(x)=
\left(
\begin{array}{c}
2\gamma(x)J\widetilde{\nu}(x) \\ 
2\left(\nabla\gamma(x)\right)^t
\end{array}
\right).
\] 
Let $\widetilde{J}\widetilde{\nu}(x)$ and 
$|J\widetilde{\nu}|(x)$be the 
cofactor matrix of the Jacobian matrix 
$J\widetilde{\nu}(x)$ and the Jacobian determinant 
of $\widetilde{\nu}$ at $x$ respectively.    
Moreover, let ${O}$ be the 
$n\times 1$ zero vector.   
Multiplying the matrix 
\[
\left(-\left(\nabla\gamma(x)\right)^t, 1\right)
\left(
\begin{array}{cc}
\widetilde{J}\widetilde{\nu}(x) &  O\\ 
O^t  & \gamma(x)|J\widetilde{\nu}|(x) 
\end{array}
\right)
\]  
to the Jacobian matrix $Jf(x)$ from the left side 
yields the following, where $E_n$ stands for the 
$n\times n$ unit matrix.   
\begin{eqnarray*}
{ } & { } & 
\left(-\left(\nabla\gamma(x)\right)^t, 1\right)
\left(
\begin{array}{cc}
\widetilde{J}\widetilde{\nu}(x) &  O\\ 
O^t  & \gamma(x)|J\widetilde{\nu}|(x) 
\end{array}
\right)
\left(
\begin{array}{c}
2\gamma(x)J\widetilde{\nu}(x) \\ 
2\left(\nabla\gamma(x)\right)^t
\end{array}
\right) \\ 
{ } & = & 
\left(-\left(\nabla\gamma(x)\right)^t, 1\right)
\left(
\begin{array}{c}
2\gamma(x)|J\widetilde{\nu}|(x) E_n \\ 
2\gamma(x)|J\widetilde{\nu}|(x)
\left(\nabla\gamma(x)\right)^t
\end{array}
\right)
\\ 
{ } & = & 
(0, \ldots, 0).
\end{eqnarray*}  
Hence we have the following.    
\begin{lemma}
Suppose that  
$|J\widetilde{\nu}|(x)\ne 0$.     
Then, we may put 
as follows: 
\begin{eqnarray*}
\nu_1(x)  & = &   
-\left(\widetilde{J}\widetilde{\nu}(x)\right)^t
\nabla \gamma(x)/||\nu(x)|| \\ 
%=
%-\left(\nabla \gamma(x)
%\widetilde{J}\widetilde{\nu}(x)
%\right)^t/||\nu(x)||
\nu_2(x)  & = &  
|J\widetilde{\nu}|(x)\gamma(x)\widetilde{\nu}(x)
/||\nu(x)||.     
\end{eqnarray*}
\end{lemma}
\noindent 
Notice that in order to show that 
$||\nu(x)||\ne 0$, the assumption  
\lq\lq $|J\widetilde{\nu}|(x)\ne 0$\rq\rq\; 
is used .   
\par 
Set {$h(x)=2g(x)-\widetilde{f}(x)$}.  
Then, 
by {elementary geometry},  
we have 
\begin{eqnarray*}
 {h(x)-P} 
 & = & f(x)-\widetilde{f}(x) \\ 
{ } & = & 
\frac{\gamma(x)}
{\widetilde{\nu}(x)\cdot \nu_2(x)}
\left(\nu_1(x)\oplus \nu_2(x)\right) \\ 
{ } & = & 
\frac{1}{|J\widetilde{\nu}|(x)}
\left(-\left(\widetilde{J}\widetilde{\nu}(x)\right)^t
\nabla \gamma(x)
\oplus |J\widetilde{\nu}|(x)\gamma(x)\widetilde{\nu}(x) \right) 
\\ 
{ } & = & 
-\frac{1}{|J\widetilde{\nu}|(x)}
\left(\widetilde{J}\widetilde{\nu}(x)\right)^t\nabla \gamma(x)
\oplus {(g(x)-P)}.
\end{eqnarray*}
Since $\widetilde{f}(x)-g(x)=g(x)-h(x)=
{(g(x)-P)}-{(h(x)-P)}$,  we have 
\[
\widetilde{f}(x)-g(x) = 
\frac{1}{|J\widetilde{\nu}|(x)}
\left(\widetilde{J}\widetilde{\nu}(x)\right)^t\nabla \gamma(x)
{\oplus {0}\widetilde{\nu}(x)}
{=}   
{\left(\left(J\widetilde{\nu}(x)\right)^{-1}\right)^t
\nabla \gamma(x)
{\oplus {0}\widetilde{\nu}(x)}}. 
\]
\hfill 
$\Box$  
%%%%%%%%%%%%%%%%%%%%%%%%%%%%%%%%%%%%%%%
\section{Application 2: Opening of Gauss mapping of 
anti-orthotomic}\label{section 5}
%%%%%%%%%%%%%%%%%%%%%%%%%%%%%%%%%%%%%%%
The application in Section \ref{section 4} is a result 
only at a non-singular point of $\widetilde{\nu}$.    
In this section, as the second application of 
Theorem \ref{theorem 1}, 
we investigate what can be asserted 
even at a singular point of $\widetilde{\nu}$.   
Several definitions are needed for the investigation of this 
section.   
\begin{definition}[\cite{ishikawaaustralia}]
\label{ramification module}
{\rm 
Let $f=(f_1, \ldots, f_n): 
(\mathbb{R}^n,0)\to (\mathbb{R}^n,0)$ be an  
equidimensional $C^\infty$ map-germ.   
\begin{enumerate}
\item[(1)] Let $\Omega_n^1$ 
denote the $\mathcal{E}_n$-module 
of $1$-forms  on $(\mathbb{R}^n,0)$.     
Then, the $\mathcal{E}_n$-module generated by $df_i$ 
$(i=1, \ldots, n)$ in $\Omega_n^1$ is called the 
\textit{Jacobi module} of 
$f$ and is denoted by {$\mathcal{J}_f$}, 
where $d h$ for 
a function-germ $h: (\mathbb{R}^n,0)\to \mathbb{R}$ 
stands for the exterior differential of $h$.    
\item[(2)]  
The \textit{ramification module} of $f$ (denoted by 
$\mathcal{R}_f$) 
is defined as the $f^*\left(\mathcal{E}_n\right)$
-module consisting of 
all function-germs $\gamma$ such that $d\gamma$ 
belongs to $\mathcal{J}_f$.   
%is an element of 
%the $\mathcal{E}_n$-module generated by 
%$df_i$ $(i=1, \ldots, n)$. 
\end{enumerate}
}
\end{definition} 
\begin{definition}[\cite{ishikawaaustralia}]
{\rm 
Let ${\color{black}\varphi} 
: (\mathbb{R}^n,0)\to (\mathbb{R}^n,0)$ be an 
equidimensional $C^\infty$ map-germ and let 
${\color{black}\delta}
: (\mathbb{R}^n,0)\to (\mathbb{R}, 0)$ a 
$C^\infty$ function-germ.    
Then, the map-germ 
$({\color{black}\varphi, \delta}):  (\mathbb{R}^n,0)\to 
(\mathbb{R}^n\times \mathbb{R},(0,0))$ 
is called an {\it opening of ${\varphi}$} if 
{$\delta\in 
\mathcal{R}_{\varphi}$}.     
}
\end{definition}
%{\color{blue}
%\begin{problem}
%{\rm 
%Let $(f, \nu): N\to \mathbb{R}^{n+1}\times S^n$ be a frontal.   
%Let $x_0\in N$ and $P\in \mathbb{R}^{n+1}$ satisfy  
%$(f(x_0)-P)\cdot \nu(x_0)\ne 0$.   
%Then, 
%is $(f-P): (N, x_0)\to \mathbb{R}^{n+1}$ 
%an opening of the Gauss map-germ 
%$\frac{f-P}{||f-P||}: (N, x_0)\to S^n$ 
%of its anti-orthotomic-germ 
%$\widetilde{f}: (N, x_0)\to \mathbb{R}^{n+1}$ ?   
%In other words, 
%is $||f-P||$ contained in  
%$\mathcal{R}_{{}_{\widetilde{\nu}}}
%=\mathcal{R}_{\left(\frac{f-P}{||f-P||}\right)}$ ? 
%}
%\end{problem}
%}
%\bigskip 
%PROBLEM 4 can be solved affirmatively as follows.   
%\end{slide}
%\begin{slide}
\begin{theorem} 
\label{theorem 3}  
%\qquad \\ 
Let $(f, \nu): N\to \mathbb{R}^{n+1}\times S^n$ be a frontal.   
Let $x_0\in N$ and $P\in \mathbb{R}^{n+1}$ satisfy  
$(f(x_0)-P)\cdot \nu(x_0)\ne 0$.   
% and 
%$\left|J\left(\frac{f-P}{||f-P||}\right)\right|(x_0)= 0$.       
Then, 
$(f-P): (N, x_0)\to \mathbb{R}^{n+1}$ 
is an opening of the Gauss map-germ   
$\widetilde{\nu}=\frac{f-P}{||f-P||}: (N, x_0)\to S^n$ 
of its anti-orthotomic     
$\widetilde{f}: (N, x_0)\to \mathbb{R}^{n+1}$.   
\end{theorem}
The optical meaning of anti-orthotomic 
is straightforward from Theorem \ref{theorem 3}.  
By definition, we have the following corollary.   
\begin{corollary}\label{corollary 5}
Let $(g, \nu): N\to \mathbb{R}^{n+1}\times S^n$ be a frontal.   
Let $x_0\in N$ and $P\in \mathbb{R}^{n+1}$ satisfy  
$(g(x_0)-P)\cdot \nu(x_0)\ne 0$.   
% and 
%$\left|J\left(\frac{f-P}{||f-P||}\right)\right|(x_0)= 0$.       
Then, 
$(g-P): (N, x_0)\to \mathbb{R}^{n+1}$ 
is an opening of the Gauss map-germ   
$\widetilde{\nu}=\frac{g-P}{||g-P||}: (N, x_0)\to S^n$ 
of its negative pedal    
$\widetilde{f}: (N, x_0)\to \mathbb{R}^{n+1}$. 
\end{corollary} 
% explains the optical meaning of 
%%%%%%%%%%%%%%%%%%%%%  
\par 
\smallskip 
\noindent 
{\bf Proof of Theorem \ref{theorem 3}.}   
\par 
Let $V\subset S^n$ be a 
sufficiently small open  neighbourhood 
of $\widetilde{\nu}(x_0)$ 
and let 
$h=(h_1, \ldots, h_n) 
: V\to T_{\widetilde{\nu}(x_0)}S^n$ be a   
normal coordinate system 
at $\widetilde{\nu}(x_0)$.
%For any $i$ $(1\le i\le n)$, set 
%\[
%\widetilde{\nu}_i=h_i\circ\widetilde{\nu}.   
%\]
Let $U$ be a sufficiently small open neighbourhood of $x_0$ 
such that $U\subset \widetilde{\nu}^{-1}(V)$ and 
set $\nu_1(x)=\left(\nu_{1,1}(x), \ldots, \nu_{1,n}(x)\right)$ 
for any $x\in U$.   
Moreover, for any $i$ $(1\le i\le n)$, set 
$\widetilde{\nu}_i=h_i\circ \widetilde{\nu}$.   
Since $\nu: N\to S^n$ is the Gauss mapping of 
$f: N\to \mathbb{R}^{n+1}$ and 
$f(x)-P=2\gamma(x)\widetilde{\nu}(x)$, we have  
\[
\sum_{i=1}^n \nu_{1, i}\gamma 
d\widetilde{\nu}_i+||\nu_2|| d\gamma=0.   
\]
Since $(f(x_0)-P)\cdot \nu(x_0)\ne 0$, it follows  
$\nu_2(x_0)\ne 0$. 
Thus, % if $U$ is sufficiently small, \\ 
we have 
\[
d\gamma = -\frac{1}{||\nu_2||}\sum_{i=1}^n 
\nu_{1,i}\gamma 
d\widetilde{\nu}_i\in \mathcal{J}_{{}_{\widetilde{\nu}}}.    
\]
\hfill 
$\Box$    
\par
\begin{figure}
\begin{center}
\includegraphics[width=8cm]
{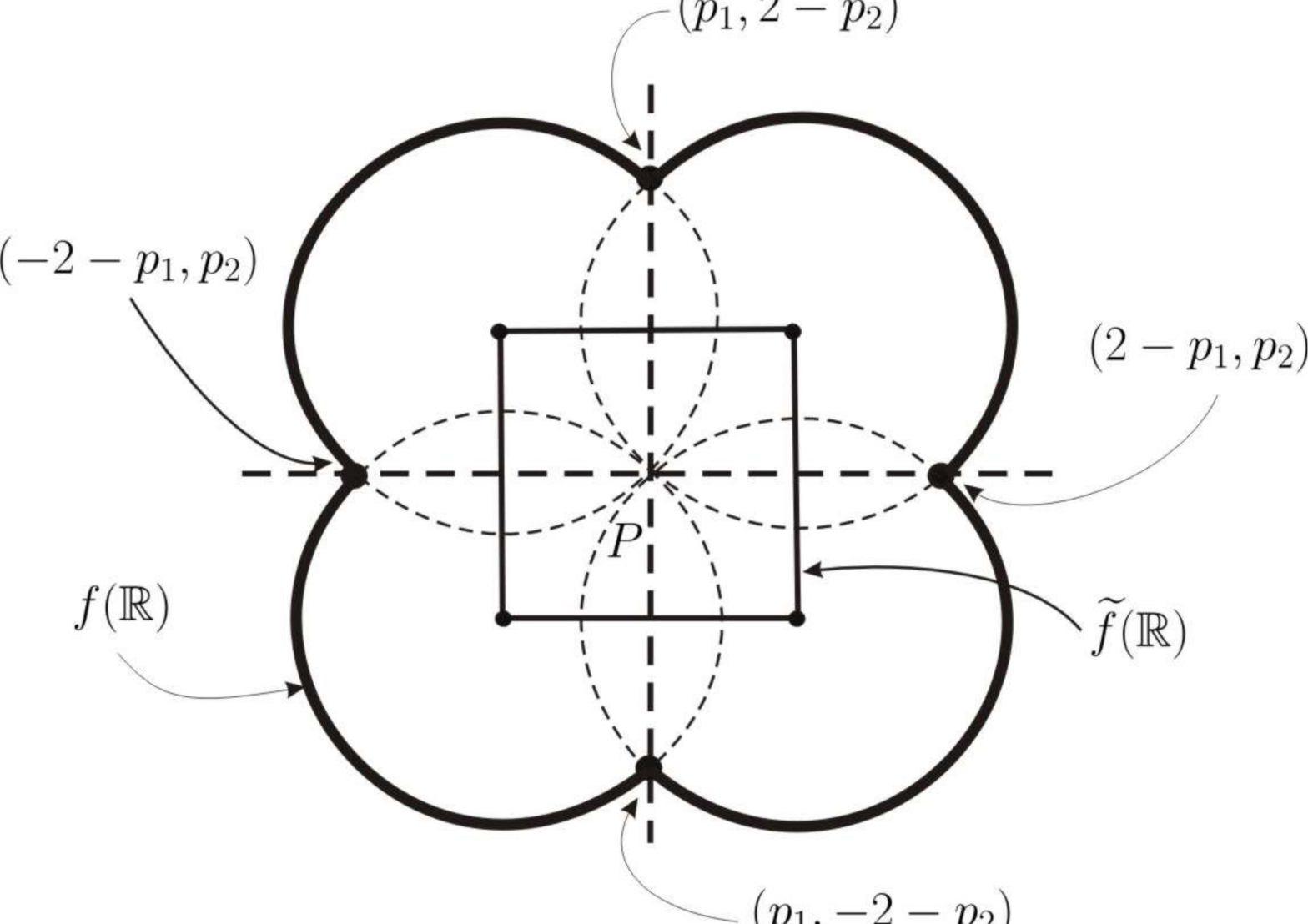}
%\qquad 
%\includegraphics[width=6.8cm]
%{orthotomic7a.jpg}
{\color{black}\caption{%Left:
How to draw $f(\mathbb{R})$ for the given square 
$\widetilde{f}(\mathbb{R})$.} 
% Right: 
}
\label{figure 5}
\end{center}
\end{figure}

 
\bigskip 
%PROPOSITION 1 and \\ 
Consider again the frontal 
$(\widetilde{f}, \widetilde{\nu}): \mathbb{R}\to 
\mathbb{R}^2\times S^1$ given in Example 
\ref{example 2}.    
Recall that the image 
$\widetilde{f}(\mathbb{R})$ is the square $S$ with 
vertexes $(1,1), (-1, 1), (1, -1), (-1, -1)$.    
Let $P=(p_1, p_2)$ be a point such that 
$-1< p_1, p_2<1$.     Then, $P$ belongs to 
$\mathcal{NS}_{\widetilde{f}}$.    
Let $f: \mathbb{R}\to \mathbb{R}^2$ be the orthotomic 
of $\widetilde{f}$ relative to $P$.    
%Although $f$ has the explicit expression,  
Theorems \ref{theorem 1} and \ref{theorem 3} 
reduce the construction of the image of $f$ 
to elementary geometry, which is explained as follows 
{\color{black}(see Figure \ref{figure 5})}.     
%%%%%%%%%%%%%%%%%%%%%%%%%%%%%%%%%%%%%%%%%%%%%%   
%%%%%%%%%%%%%%%%%%%%%%%%%%%%%%%%%%%%%%%%%%%%%%%
%%%%%%%%%%%%%%%%%%%%%%%%%%%%%%%%%%%%%%%%%%% 
%%%%%%%%%%%%%%%%%%%%%%%%%%%%%%%%%%%%%%%%%%% 
%\vspace{7cm}
%%%%%%%%%%%%%%%%%%%%%%%%%%%%%%%%%%%%%%%%%%% 
%%%%%%%%%%%%%%%%%%%%%%%%%%%%%%%%%%%%%%%%%%%  
By the construction of $\widetilde{\nu}$, 
if $t$ belongs to the union of closed intervals 
\[
\bigcup_{n\in \mathbb{Z}} \left([8n+1, 8n+2]
\cup [8n+3, 8n+4]\cup [8n+5, 8n+6]\cup 
[8n+7, 8n+8]\right), 
\] 
then $t$ must be a singular point of $\widetilde{\nu}$.   
Thus, by Theorem \ref{theorem 3}, $t$ 
must be a singular point of $f$, and therefore 
each of $f([8n+1, 8n+2]), f([8n+3, 8n+4]), f([8n+5, 8n+6]), 
f([8n+7, 8n+8])$ must be one point.     
By definition, the one point must be the mirror 
image of $P$ as follows.    
\begin{lemma}
\begin{eqnarray*}
f([8n+1, 8n+2])=(2-p_1, p_2), & f([8n+3, 8n+4])=(p_1, 2-p_2), \\ 
f([8n+5, 8n+6])=(-2-p_1, p_2), & 
f([8n+7, 8n+8])=(p_1, -2-p_2).   
\end{eqnarray*}
\end{lemma} 
By the construction of $\widetilde{f}$, 
\begin{eqnarray*}
\widetilde{f}([8n, 8n+1])=(1,-1), 
& \widetilde{f}([8n+2, 8n+3])=(1,1), \\ 
\widetilde{f}([8n+4, 8n+5])=(-1, 1), & 
\widetilde{f}([8n+6, 8n+7])=(-1, -1).   
\end{eqnarray*} 
By the assertion (4) of Theorem \ref{theorem 1}, 
the following holds.   
\begin{lemma}
\[
||\widetilde{f}(t)-f(t)|| = 
\left\{
\begin{array}{cc}
\sqrt{(p_1-1)^2+(p_2+1)^2} & (\mbox{if }8n\le t\le 8n+1), \\ 
\sqrt{(p_1-1)^2+(p_2-1)^2} & (\mbox{if }8n+2\le t\le 8n+3), \\
\sqrt{(p_1+1)^2+(p_2-1)^2} & (\mbox{if }8n+4\le t\le 8n+5), \\
\sqrt{(p_1+1)^2+(p_2+1)^2} & (\mbox{if }8n+6\le t\le 8n+7). \\
\end{array}
\right.
\]
\end{lemma}
By the construction of $\widetilde{\nu}$, 
we have the following.  
\begin{lemma}
\begin{enumerate}
\item[(1)]\quad 
$f([8n, 8n+1])$ is exactly the hemicircle centered at 
$\widetilde{f}([8n, 8n+1])=(1, -1)$ with boundary 
$f([8n+7, 8n+8])=(p_1, -2-p_2)$ and 
$f([8n+1, 8n+2])=(2-p_1, p_2)$ which does not 
contain $P$.   
\item[(2)]\quad 
$f([8n+2, 8n+3])$ is exactly the hemicircle centered at 
$\widetilde{f}([8n+2, 8n+3])=(1, 1)$ with boundary 
$f([8n+1, 8n+2])=(2-p_1, p_2)$ 
and $f([8n+3, 8n+4])=(p_1, 2-p_2)$ which does not contain 
$P$.  
\item[(3)]\quad 
$f([8n+4, 8n+5])$ is exactly the hemicircle centered at 
$\widetilde{f}([8n+4, 8n+5])=(-1, 1)$ with boundary 
$f([8n+3, 8n+4])=(p_1, 2-p_2)$ 
and $f([8n+5, 8n+6])=(-2-p_1, p_2)$ which does not 
contain $P$.   
\item[(4)]\quad 
$f([8n+6, 8n+7])$ is exactly the hemicircle centered at 
$\widetilde{f}([8n+6, 8n+7])=(-1, -1)$ with boundary 
$f([8n+5, 8n+6])=(-2-p_1, p_2)$ 
and $f([8n+7, 8n+8])=(p_1, -2-p_2)$ which does not 
contain $P$.   
\end{enumerate}
\end{lemma}
For the precise shape of the pedal $g: \mathbb{R}\to 
\mathbb{R}^2$ of $\widetilde{f}$ relative to $P$, 
just shrink $f(\mathbb{R})$ to $50$ percent 
with respect to  $P$.      
\par 
It seems that the method of C.~Herring explained 
in \cite{herring} is similar as our method.   
However, his method seems to rely 
on a thermodynamical consideration of atoms.     
Our method needs no physical consideration.   
Once the given shape is realized 
as the image of frontal $\widetilde{f}(\mathbb{R})$, 
by applying Theorem \ref{theorem 1} 
and Theorem \ref{theorem 3}, only elementary geometry 
is needed.    
In other words, under any physical situation, 
if the same square is given, 
then the $\gamma$-plot for the square 
must have the same shape.   
%only elementary geometry is needed 
%for the construction of precise shape of $f(\mathbb{R})$.   
%%%%%%%%%%%%%%%%%%%%%%%%%%%%%%%%%%%%%%%
%%%%%%%%%%%%%%%%%%%%%%%%%%%%%%%%%%%%%%% 
\section{Application 3: A criterion for fronts}\label{section 6}
%%%%%%%%%%%%%%%%%%%%%%%%%%%%%%%%%%%%%%% 
%%%%%%%%%%%%%%%%%%%%%%%%%%%%%%%%%%%%%%% 
\begin{definition}
{\rm A germ of frontal  
$(f, \nu): (N, x_0)\to \mathbb{R}^{n+1}\times S^n$ 
is called a {\it germ of front} (or {\it front-germ}) if 
$(f, \nu)$ is non-singular at $x_0$.   
}
\end{definition} 
\noindent 
Given a frontal 
$(f, \nu): N\to \mathbb{R}^{n+1}\times S^n$, 
if 
$(f, \nu): (N, x)\to \mathbb{R}^{n+1}\times S^n$ is 
a germ of front for any $x\in N$, then 
$(f, \nu)$ is called a {\it front}.    
A front is also called a {\it wave-front}.   
For details on fronts, 
see for example \cite{arnold, arnoldetall}.  
%This problem can be solved as follows.   
\begin{theorem}\label{theorem 4}
Let $(f, \nu): N\to \mathbb{R}^{n+1}\times S^n$ be a frontal 
and let $x_0$ be a point of $N$.     
Then, for any point $P\in\mathbb{R}^{n+1}$ such that 
$(f(x_0)-P)\cdot \nu(x_0)\ne 0$, the following are equivalent, 
where $\left(\widetilde{f}, \widetilde{\nu}\right): 
(N, x_0)\to \mathbb{R}^{n+1}\times S^n$ 
is the anti-orthotomic germ of $f$ relative to $P$.    
\begin{enumerate}
\item[(1)] $(f, \nu): (N,x_0)\to \mathbb{R}^{n+1}\times S^n$ 
is a front-germ.   
\item[(2)] 
$\left(\widetilde{f}, \widetilde{\nu}\right): 
(N,x_0)\to \mathbb{R}^{n+1}\times S^n$ 
is a front-germ.  
\item[(3)] 
$\left(f, \widetilde{f}\right): (N,x_0)\to 
\mathbb{R}^{n+1}\times \mathbb{R}^{n+1}$ is non-singular.   
\end{enumerate}   
\end{theorem}
Theorem \ref{theorem 4} answers the question 
communicated by 
A.~Honda and K.~Teramoto (\cite{hondateramoto}).     
%By Proposition \ref{proposition 1},   
Theorem \ref{theorem 4} 
yields the following corollaries.   
%, Theorem \ref{theorem 3}, 
%Theorem \ref{theorem 4}, 
%Corollary \ref{corollary 1, corollary 2, corollary 5}, 
%we have the following  corollaries.   
\begin{corollary}
Let $(g, \nu): N\to \mathbb{R}^{n+1}\times S^n$ be a frontal
and let $x_0$ be a point of $N$.     
Then, for any point $P\in\mathbb{R}^{n+1}$ such that 
$(g(x_0)-P)\cdot \nu(x_0)\ne 0$, the following are equivalent,  
where $\left(\widetilde{f}, \widetilde{\nu}\right): 
(N, x_0)\to \mathbb{R}^{n+1}\times S^n$ 
is the negative pedal germ of $g$ relative to $P$.  
\begin{enumerate}
\item[(1)] $(g, \nu): (N,x_0)\to \mathbb{R}^{n+1}\times S^n$ 
is a front-germ.   
\item[(2)] 
$\left(\widetilde{f}, \widetilde{\nu}\right): 
(N,x_0)\to \mathbb{R}^{n+1}\times S^n$ 
is a front-germ.  
\item[(3)] 
$\left(g, \widetilde{f}\right): (N,x_0)\to 
\mathbb{R}^{n+1}\times \mathbb{R}^{n+1}$ is non-singular.   
\end{enumerate}   
\end{corollary}
\begin{corollary}
Let $(f, \nu): N\to \mathbb{R}^{n+1}\times S^n$ be a frontal.    
Let two points $x_0\in N$ and 
$P\in\mathbb{R}^{n+1}$ satisfy   
$(f(x_0)-P)\cdot \nu(x_0)\ne 0$.    
Let $(\widetilde{f}, \widetilde{\nu}): 
N\to \mathbb{R}^{n+1}\times S^n$ 
be the anti-orthotomic of $f$ relative to $P$.   
If $x_0$ is not contained in 
Sing$(\widetilde{f})\cap$Sing$(\widetilde{\nu})$, 
then the map-germ $f: (N, x_0)\to \mathbb{R}^{n+1}$ 
is a front-germ; 
where for a $C^\infty$ mapping $\varphi: X\to Y$, 
Sing$(\varphi)$ stands for the singular set 
of $\varphi$.     
In particular, if %at least 
%one of 
$\widetilde{f}: 
(N, x_0)\to \mathbb{R}^{n+1}$ 
%and $\widetilde{\nu}: N\to S^n$ 
is non-singular, then $f: (N, x_0)\to \mathbb{R}^{n+1}$ 
must be a front-germ.       
\end{corollary}
\begin{corollary}
Let $(g, \nu): N\to \mathbb{R}^{n+1}\times S^n$ be a frontal.    
Let two points $x_0\in N$ and 
$P\in\mathbb{R}^{n+1}$ satisfy   
$(g(x_0)-P)\cdot \nu(x_0)\ne 0$.    
Let $(\widetilde{f}, \widetilde{\nu}): 
N\to \mathbb{R}^{n+1}\times S^n$ 
be the negative pedal of $g$ relative to $P$.   
If $x_0$ is not contained in 
Sing$(\widetilde{f})\cap$Sing$(\widetilde{\nu})$, then     
the map-germ $g: (N, x_0)\to \mathbb{R}^{n+1}$ 
is a front-germ.     
In particular, if %at least 
%one of 
$\widetilde{f}: 
(N, x_0)\to \mathbb{R}^{n+1}$ 
%and $\widetilde{\nu}: N\to S^n$ 
is non-singular, then $g: (N, x_0)\to \mathbb{R}^{n+1}$ 
must be a front-germ.       
\end{corollary}
\begin{corollary}
Let $(\widetilde{f}, \widetilde{\nu}): 
N\to \mathbb{R}^{n+1}\times S^n$ be a frontal.    
Let two points $x_0\in N$ and 
$P\in\mathbb{R}^{n+1}$ satisfy   
$(\widetilde{f}(x_0)-P)\cdot \widetilde{\nu}(x_0)\ne 0$.    
Let $({f}, {\nu}): 
N\to \mathbb{R}^{n+1}\times S^n$ 
be the orthotomic of $\widetilde{f}$ relative to $P$.   
If $x_0$ is not contained in 
Sing$({f})\cap$Sing$({\nu})$, then     
the map-germ $\widetilde{f}: (N, x_0)\to \mathbb{R}^{n+1}$ 
is a front-germ.         
%In particular, if at least 
%one of 
%${f}: (N, x_0)\to \mathbb{R}^{n+1}$ 
%and ${\nu}: (N, x_0)\to S^n$ 
%is non-singular, 
%then $\widetilde{f}: (N, x_0)\to \mathbb{R}^{n+1}$ 
%must be a front-germ.       
\end{corollary}
\begin{corollary}
Let $(\widetilde{f}, \widetilde{\nu}): 
N\to \mathbb{R}^{n+1}\times S^n$ be a frontal.    
Let two points $x_0\in N$ and 
$P\in\mathbb{R}^{n+1}$ satisfy   
$(\widetilde{f}(x_0)-P)\cdot \widetilde{\nu}(x_0)\ne 0$.    
Let $({g}, {\nu}): 
N\to \mathbb{R}^{n+1}\times S^n$ 
be the pedal of $\widetilde{f}$ relative to $P$.   
If $x_0$ is not contained in 
Sing$({g})\cap$Sing$({\nu})$, then     
the map-germ $\widetilde{f}: (N, x_0)\to \mathbb{R}^{n+1}$ 
is a front-germ.           
%In particular, if at least 
%one of 
%${g}: (N, x_0)\to \mathbb{R}^{n+1}$ 
%and ${\nu}: (N, x_0)\to S^n$ 
%is non-singular, 
%then $\widetilde{f}: (N, x_0)\to \mathbb{R}^{n+1}$ 
%must be a front-germ.       
\end{corollary}
\par 
\medskip 
\noindent 
{\bf Proof of Theorem \ref{theorem 4}.} 
\par  
%\proof\quad 
%{\bf Proof of Theorem 4}. \\  
$(2)\Leftrightarrow (3)$ is just a corollary of Theorem  
\ref{theorem 3}.  
%By Corollary 3, it follows that $(2)\Leftrightarrow (3)$.   
$(1)\Rightarrow (3)$ is trivial.   
% and $(2)\Rightarrow (3)$.  
%In order to show $(3)\Rightarrow (1)$, 
%the no-silhouette condition is used.   
%$(3)\Rightarrow (1)$ is proved as follows.   \\ 
%\smallskip 
Thus, in order to complete the proof, it is 
sufficient to show $(3)\Rightarrow (1)$.   
% can be proved in the following way.   \\ 
Suppose that $(f, \widetilde{f})$ is non-singular.   
Then, $(f, f-\widetilde{f})$ is non-singular.    
By the assumption 
$({f}(x_0)-P)\cdot {\nu}(x_0)\ne 0$,  
%by Proposition \ref{proposition 1}, 
%$P\in \mathcal{NS}_f$, 
%By the no-silhouette condition, 
it follows that the projection  
$
\pi: d_{x_0}(f-\widetilde{f})(\mbox{\rm Ker}(d_{x_0} f))
\to T_{{\nu(x_0)}}S^n
$  
is injective.   
Therefore, $\left(f, {\nu}\right)=
\left(f, \frac{f-\widetilde{f}}{||f-\widetilde{f}||}\right)
: (N, x_0)\to \mathbb{R}^{n+1}\times S^n$ 
is non-singular.   %\\ 
%$(3)\Rightarrow (2)$ can be proved similarly.   
\hfill $\Box$
%%%%%%%%%%%%%%%%%%%%%%%%%%%%%%%%%%%%%%%
%%%%%%%%%%%%%
%%%%%%%%%%%%%% %%%%%%%%%%%%%%%%%%%%%%%%%%%%%%%%%%%%%
\section*{Acknowledgement}
%\begin{acknowledgements}
%\thanks
The authors thank Masatomo Takahashi for teaching them 
his improvement of classical envelope.   
The second author was %partially 
supported
by JSPS KAKENHI Grant Number 17K05245.
%\end{acknowledgements}
%%%%%%%%%%%%%%%%%
%%%%%%%%%%%%%%%%%%%%%%%%%%%%%%%%%%%

\end{document}